\documentclass[11pt,a4paper]{article}
\usepackage[utf8]{inputenc}
\usepackage{amssymb, amsmath, mathtools}
\usepackage[top = 2.5cm, bottom = 2.5cm, left = 2.5cm, right = 2.5cm]{geometry} 
\usepackage{enumitem}
\usepackage{commath}
\usepackage{color}
\usepackage{mathrsfs}
\usepackage{graphicx}
\usepackage{svg}
\usepackage{dsfont}
\usepackage{systeme}
\usepackage{amsmath}
\usepackage{hyperref}
\usepackage{bbold}
\usepackage{comment}
\usepackage{subcaption}
\usepackage{float}
\usepackage{cleveref}
\usepackage[toc,page]{appendix}
\usepackage[dvipsnames]{xcolor}
\usepackage{authblk}

\usepackage[
  backend=biber,
  style=numeric-comp,
  sorting=none,
  maxbibnames=99
]{biblatex}
\author[1,2, *]{Aurora Faure Ragani}
\author[1,2,3]{Robbin Bastiaansen}

\affil[1]{Mathematical Institute, Utrecht University, Utrecht, the Netherlands}
\affil[2]{Institute for Marine and Atmospheric Research Utrecht (IMAU), Utrecht University, Utrecht, the Netherlands}
\affil[3]{Department of Physics, Centre for Complex System Studies, Utrecht University, Utrecht, The Netherlands}
\affil[*]{Corresponding author: \href{mailto:a.faureragani@uu.nl}{a.faureragani@uu.nl}}

\begin{document}
\title{Optimisation of tipping pathways in a spatially heterogeneous world}
\maketitle

\begin{abstract}
In spatially extended systems, tipping does not necessarily lead to a uniform, abrupt transition typical of low-dimensional conceptual climate models. Instead, spatially structured forms of tipping can emerge, for example, through front propagation and pinning, coexistence between alternative states, and pattern formation. Tipping can therefore remain partial or spatially heterogeneous rather than affecting the entire system.
Existing work on spatial tipping in conceptual models remains mostly descriptive, and a general framework to influence tipping dynamics through spatial interventions is still lacking.
Here, we introduce a constrained optimisation framework that systematically identifies spatial interventions designed to maximise resilience or promote recovery of a desirable state, subject to dynamical and resource constraints. The framework is general and applicable to a broad class of spatially extended systems.
We illustrate it using the one-dimensional Allen–Cahn equation with spatial heterogeneity, a minimal bistable model in which the key spatial tipping processes in gradient systems can be analysed explicitly. In this setting, optimisation acts by shifting local tipping thresholds and controlling front dynamics.
Our results show that small-scale local, targeted modifications can determine large-scale system outcomes: optimal interventions can prevent tipping, induce recovery of the desirable state, or confine collapse to limited regions of the domain.
\end{abstract}

\section{Introduction}
In the last 20 years, many Earth subsystems have been classified as tipping systems \cite{lenton2008_tipping_elements, IPCC2023, lenton2019_climate_tipping_points}. These systems can respond nonlinearly, and sometimes abruptly, to gradual changes in external forcing. Such abrupt responses constitute tipping events, in which the system shifts to an alternative state that is often less desirable \cite{scheffer2001_catastrophic_shifts, armstrongMcKay2022_tipping_points}. Once in this undesirable state, the system may remain there even if external conditions improve. This behaviour, known as hysteresis, prevents easy recovery to the desirable state. Because such irreversible shifts can profoundly alter Earth systems, it is crucial to understand, predict, and manage these transitions.

In conceptual studies of these transitions, low-dimensional models of tipping systems are often used; these models usually give a binary view of tipping: the system is either in a desirable state or in an undesirable one. In such models, tipping takes the form of a single abrupt shift in which the entire system rearranges. Real systems, however, differ from these simplified models because they possess many more degrees of freedom. A major source of this additional complexity is the spatial dimension \cite{meron2015_nonlinear_ecosystems, banerjee2026_rethinking_spatial_tipping}: when spatial structure is included, more subtle forms of tipping can emerge \cite{bel2012_gradual_regime_shift_spatially_ext_syst, zelnik2018_regime_shift_front_dynamics, champneys2021_pinning}. These include, for instance, coexistence states \cite{siteur2014_bistability_spatial_resilience}, fragmented tipping pathways \cite{bastiaansen2022_fragmented_tipping}, and pattern-mediated transitions \cite{rietkerk2008_patterns_in_ecosystems, rietkerk2021_evasion_tipping}. Because spatial structure reshapes how systems tip, it also raises new questions about how these transitions might be influenced in practice.

Given the potentially severe consequences of tipping events, a central question is whether they can be prevented or mitigated. Most existing work, however, has focused mainly on understanding tipping behaviour rather than shaping it. Conceptual models have helped identify the mechanisms underlying abrupt transitions, often revealing bistability and saddle-node bifurcations that separate alternative stable states \cite{scheffer2001_catastrophic_shifts, lenton2008_tipping_elements}. 

Only a few studies have explored intervention strategies in tipping-prone systems \cite{zelnik2021_human_intervention_in_ecosystems, vidiella2018_delayed_transition_exploitation}, but a general framework remains lacking, particularly for spatially extended systems where interventions must be targeted in space and implemented under resource constraints. As a result, fundamental questions remain open: how can we maximise a system resilience, and how can we determine the most effective spatial intervention under certain constraints given by society's capabilities?

Many real-world systems share the spatial structure and multistability that motivate our approach \cite{bastiaansen2022_fragmented_tipping}. Dryland ecosystems provide a clear example, as they often exhibit fronts separating vegetated and bare states, with degradation or recovery spreading spatially \cite{rietkerk2002_selforganized_vegetation,meron2016_pattern_missing_link,meron2018_dryland_patterns}. Similar spatially extended multistable dynamics have also been discussed for other systems, including the Atlantic Meridional Overturning Circulation (AMOC) \cite{lohmann2024_multistability_amoc}, grounding lines in marine ice sheets \cite{pegler2018_ice_sheet}, and eutrophication fronts in shallow lakes or coastal waters \cite{scheffer2001_catastrophic_shifts}. In such settings, the possibility of targeted spatial interventions motivates the need for a general optimisation framework.

In this paper, we address this gap and develop an optimisation framework that identifies spatial interventions designed to maximise resilience and enhance recovery of a tipping-prone system. In this framework, an intervention is represented as a spatial perturbation of the system parameters. The effect of such interventions is evaluated using a resilience measure. This measure is defined by a scalar objective function that quantifies the system's ability to maintain or recover a desirable state. The precise form of this function depends on the system and the question under consideration. The flexibility in choosing the objective function reflects the multiple ways to quantify resilience in dynamical systems, depending on the system itself, the perturbation considered, and the aspect of persistence or recovery of interest \cite{krakovska2024_resilience_dynamical_systems, holling1973_resilience, scheffer2001_catastrophic_shifts}

Within this framework, we consider two complementary optimisation strategies. The first, which we call equilibrium optimisation, targets the bifurcation structure, revealed by continuation, to delay a tipping point. The second, which we call final-state optimisation, steers time-dependent trajectories to halt front propagation, promote recovery, or partially prevent tipping when complete avoidance is not possible. We formulate both strategies as constrained optimisation problems in which an objective function quantifies resilience and some constraints reflect limited resources and impose the governing dynamics. Although we illustrate the framework using gradient systems—because they provide explicit tipping thresholds and well-understood front dynamics—the framework is general and applicable to a broader class of spatial models.

Many spatial tipping systems arise as spatial extensions of the low-dimensional conceptual models discussed above. Ignoring spatial effects, the local dynamics of a state variable $y$ can be described by an ordinary differential equation.
When spatial interactions are included, diffusion is commonly used to represent local coupling between neighbouring locations. This leads naturally from a conceptual tipping model to its spatial extension:
\[
\frac{dy}{dt}=f(y;\mu)
\;\;\longrightarrow\;\;
\frac{\partial y}{\partial t}=D\frac{\partial^2 y}{\partial x^2}+f(y;\mu)
\]
Such reaction–diffusion equations capture key spatial tipping behaviour \cite{meron2015_nonlinear_ecosystems, zelnik2018_regime_shift_front_dynamics, bastiaansen2022_fragmented_tipping}. In this paper, we focus on spatially extended gradient systems, a tractable class of reaction–diffusion models whose dynamics can be expressed in terms of a potential landscape, allowing tipping thresholds and front dynamics to be characterised explicitly.

The rest of the paper is structured as follows.
\Cref{sec:2} introduces the theory of gradient systems, with a focus on front solutions and their relation to the Maxwell point and tipping thresholds.
\Cref{sec:3} presents the optimisation framework in detail and describes the two techniques we employ: equilibrium optimisation and final-state optimisation.
\Cref{sec:4} applies the framework to the Allen–Cahn equation, a simple bistable PDE that captures the key mechanisms of gradient systems. Three examples are analysed there: one using equilibrium optimisation and two using final-state optimisation.
Finally, \Cref{sec:5} discusses implications, limitations, and generality of the framework.

\section{Theory}\label{sec:2}
We consider a spatially extended gradient system with a state variable $y(x, t)$ on a one-dimensional spatial domain:
\begin{equation}\label{eq:AC}
\frac{\partial y}{\partial t} = D \frac{\partial^2 y}{\partial x^2} + f(y; \mu) = D \frac{\partial^2 y}{\partial x^2} - \frac{\partial V}{\partial y}(y, \mu),
\end{equation}
where $D$ is the diffusion coefficient, $\mu$ the bifurcation parameter, and $V(y, \mu)$ the potential associated with the local dynamics.
In this section, we focus on those properties of spatially extended gradient systems that determine tipping and front behaviour, namely the location of saddle-node bifurcations and the motion and pinning of fronts at the Maxwell point.

First, \eqref{eq:AC} admits spatially homogeneous equilibria, corresponding to solutions that are constant in space. For fixed $\mu$, these equilibria satisfy $f(y;\mu) = 0$, and therefore correspond to the critical points of the potential $V$ $\bigl( \frac{\partial V}{\partial y} = 0\bigr)$. 

\begin{figure}
    \centering
        \includegraphics[width=0.9\linewidth]{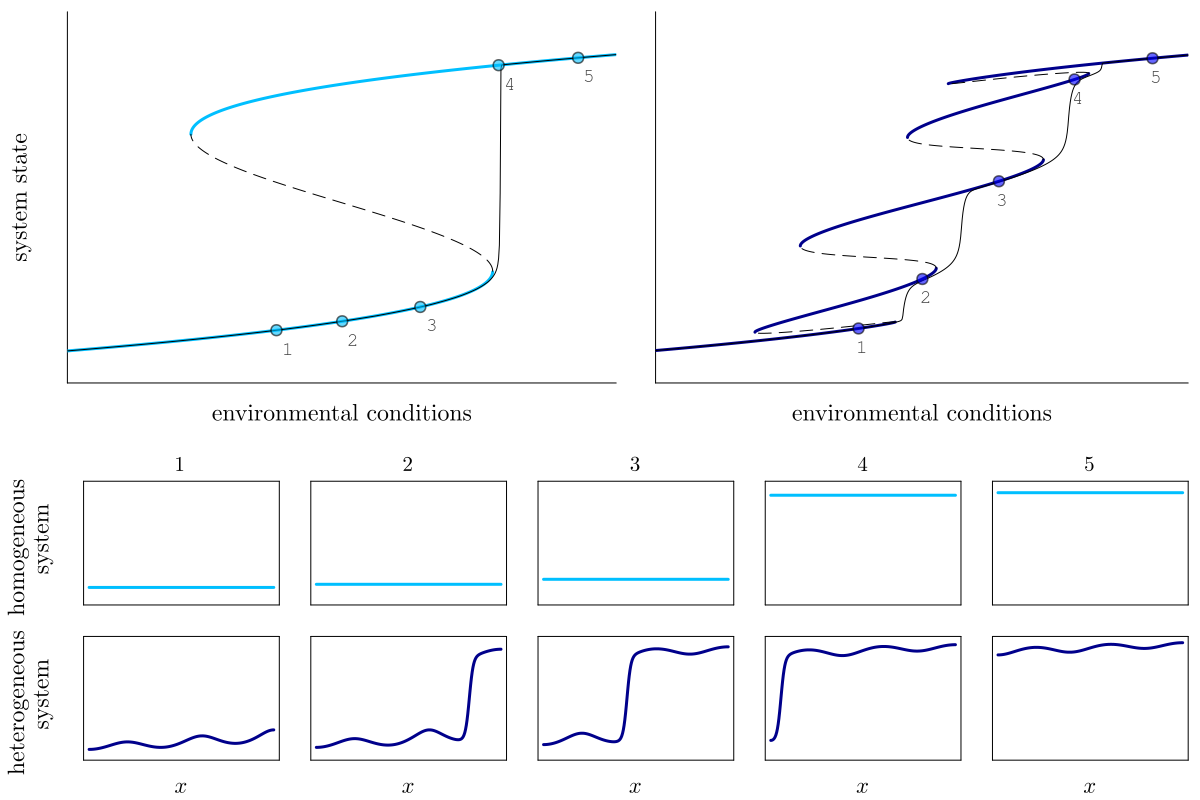}
    \caption{Representative bifurcation diagrams of the system \eqref{eq:AC} (left) and \eqref{eq:AC2} (right). Thick coloured curves indicate stable equilibrium branches, and dashed curves indicate unstable branches. The thin black lines show the system trajectory along the bifurcation structure, starting from the lower branch. Steady-state profiles, corresponding to the dots along the branches, are shown below the bifurcation diagrams.}\label{fig:0}
\end{figure}

In addition to these homogeneous equilibria, the system also admits spatially heterogeneous equilibrium states, often called coexistence states, in which different regions of the domain lie in different homogeneous equilibria. These regions are connected by spatial interfaces, or fronts, whose motion and stability determine the large-scale dynamics of the system \cite{banerjee2026_rethinking_spatial_tipping, CarrPego1989_PDEpatterns}(see the bottom panels in \Cref{fig:0} for examples).

To analyse front motion explicitly, we restrict attention to the simplest bistable case in which the potential $V(y, \mu)$ has two minima $y^-$ and $y^+$, for a certain range of the parameter $\mu$.
A travelling front solution connecting these two states can be written as $y(x,t) = u(x-ct)$, where $c$ denotes the front velocity.
To determine $c$ explicitly, we substitute this travelling-wave form into the governing equation \eqref{eq:AC}, reducing the PDE to the following ODE:
\[
-c\,\frac{du}{dz} = D\frac{d^2 u}{dz^2} + f(u;\mu).
\]
By multiplying by $\frac{du}{dz}$ and integrating the above equation against $dz$ over the real line,
one obtains the front speed:
\[
c = \frac{V(y^-;\mu)-V(y^+;\mu)}
{\displaystyle \int_{-\infty}^{\infty}\left(\frac{du}{dz}\right)^2 dz}.
\]
Because the denominator is strictly positive for a nontrivial front, the sign of $c$ is determined by the potential difference $V(y^-;\mu)-V(y^+;\mu)$. The parameter value $ \mu_M$ at which the two minima have equal depth ($V(y^-;\mu)=V(y^+;\mu)$) is called the Maxwell point. At this point, neither state invades the other, and the interface remains stationary ($c = 0$); in other words, the front is pinned. This leads to an additional equilibrium of \eqref{eq:AC}. 

Front motion depends on differences in the local potential. For this reason, spatial variation in local conditions can significantly change the behaviour of the system. In many real systems, environmental factors such as soil properties, topography, or resource availability vary across space. These variations lead to location-dependent dynamics.
To represent this spatial variation in front dynamics, we allow the reaction term of \eqref{eq:AC} to depend explicitly on position,
\begin{equation}\label{eq:AC2}
\frac{\partial y}{\partial t}
= D \frac{\partial^2 y}{\partial x^2} + f(y, x; \mu)
= D \frac{\partial^2 y}{\partial x^2} - \frac{\partial V}{\partial y}(y, x, \mu).
\end{equation}
With this formulation, key parameter values (tipping points and Maxwell point) can differ across space, so that different regions of the domain may favour different homogeneous states even when $\mu$ is fixed.
Moreover, front dynamics becomes local, with both speed and direction of propagation depending on the potential differences at the front location.
When different regions favour different states, coexistence states can involve several fronts. If such fronts are close, they may interact with each other and eventually annihilate each other \cite{banerjee2026_rethinking_spatial_tipping, bastiaansen2025_multifrontdynamics_spatially_inhomogeneous}.
All of this motivates the implementation of spatial modifications that locally shift the bifurcation parameter to target tipping thresholds and front pinning.

The effects of spatial heterogeneity become apparent when the system dynamics are represented in a bifurcation diagram \cite{bastiaansen2022_fragmented_tipping}. \Cref{fig:0} compares the bifurcation structure of a bistable gradient system without (left panel) and with (right panel) spatial heterogeneity. In the homogeneous case, the diagram exhibits the familiar S-shaped structure: two stable branches separated by an unstable branch, with saddle–node bifurcations marking the boundaries of the bistable regime.
In the heterogeneous system, the bifurcation parameter varies across space. The resulting bifurcation diagram contains additional stable branches corresponding to stable coexistence states. Along these branches, different regions of the domain lie near different homogeneous equilibria and are separated by stationary fronts. Snapshots of the system state along the branches illustrate these configurations.
The presence of these intermediate branches changes the tipping behaviour of the system. Crossing a tipping point does not necessarily trigger a transition of the entire domain. Instead, the system may move to a coexistence state in which only part of the domain shifts to the alternative equilibrium. Tipping, therefore, occurs in a partial and fragmented manner.
The precise structure of the bifurcation diagram depends on the spatial shape and magnitude of the heterogeneity \cite{bastiaansen2022_fragmented_tipping}. Nevertheless, the qualitative mechanisms described above are general and apply to any spatially heterogeneous gradient system. Because spatial heterogeneity alters both tipping thresholds and front dynamics, appropriately designed spatial perturbations can influence system resilience. This motivates the optimisation framework introduced in the next section.

\section{Optimisation Framework}\label{sec:3}
We formulate resilience optimisation as a constrained mathematical optimisation problem. In this formulation, the decision variable $z(x)$ represents a spatially heterogeneous intervention on the system's parameters, and its admissible forms define the search space $\Omega$. The objective is to identify an intervention $z \in \Omega$ that optimises a scalar objective function $F(z)$, representing system resilience, while satisfying constraints that encode both the system dynamics and practical limitations on the intervention. In general, the optimisation problem can be written as
\begin{align}
    &\min_{z \in \Omega} F(z) \label{eq:opt_general} \\
    &\text{subject to } 
    \begin{cases}
        G_i(z) \le 0, & i = 1, \dots, m, \\
        H_j(z) = 0,   & j = 1, \dots, p.
    \end{cases} \nonumber
\end{align}
Here, $F$ defines the optimisation objective, while $G_i$ and $H_j$ denote inequality and equality constraints. Inequality constraints typically capture practical limitations, such as finite resources or physical feasibility, whereas equality constraints enforce the governing equations.  

Concretely, the objective function specifies which aspect of resilience the intervention is designed to improve. 
Depending on the question, \(F\) may quantify a distance to a tipping threshold, a distance to the boundary of the basin of attraction of the desirable state, the return time after a finite perturbation, or the expected escape time under stochastic forcing. It may also be defined directly from a model trajectory, for example, by measuring the extent of the spatial domain that remains in, or recovers to, the desirable state over a prescribed time horizon. 
Thus, \(F\) translates a qualitative intervention goal, such as delaying tipping, increasing tolerance to disturbances, or reducing the spatial extent of collapse, into a scalar quantity that can be optimised subject to the governing dynamics and intervention constraints.

\begin{figure}
    \centering
    \includegraphics[width=\linewidth]{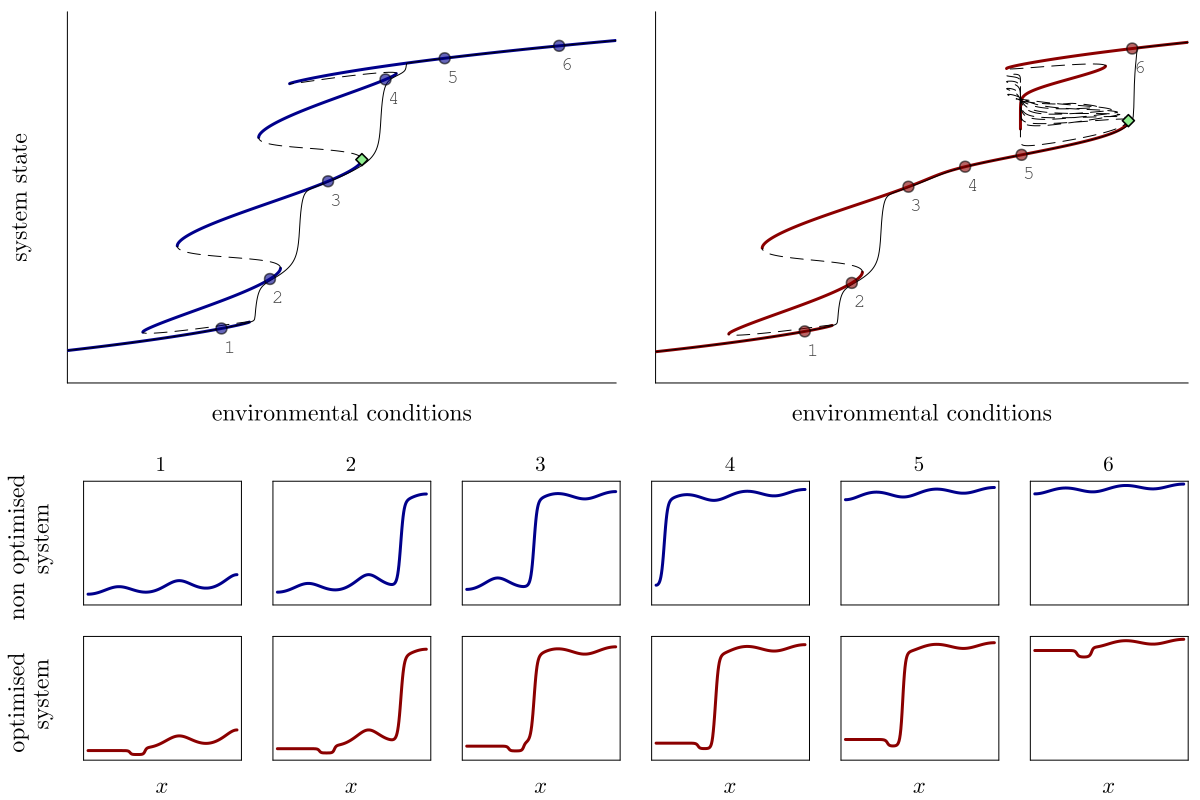}
    \caption{Bifurcation diagrams before (left) and after optimisation (right), together with representative spatial snapshots of the corresponding trajectories. The top row shows the original system (left) and the optimised system (right). The horizontal axis denotes environmental conditions, with larger values corresponding to more adverse conditions, while the vertical axis represents the system state. Thick coloured curves indicate stable equilibrium branches, and dashed curves indicate unstable branches. The thin black curve shows the trajectory followed by the system along the bifurcation structure, starting from the lower branch. Dots mark the environmental conditions at which the snapshots shown in the lower panels are taken; the same values are used in both diagrams. In the non-optimised case, successive saddle-node bifurcations separate the marked states. The first transition (from state 1 to state 2) corresponds to the formation of a front. Subsequent transitions (from state 2 to state 3 and from state 3 to state 4) are associated with abrupt leftward shifts of the front, while the transition from state 4 to state 5 corresponds to the disappearance of the front. The green diamond marks the saddle node selected as the optimisation target, with the aim of preventing the transition from state 3 to state 4. In the optimised case, a spatial perturbation shifts the targeted saddle node towards more adverse environmental conditions. As a result, the transition from state 3 to state 4 is avoided over the range shown, and the subsequent transition from state 4 to state 5 is likewise avoided.}
    \label{fig:FIGURE1}
\end{figure}

Once this objective has been specified, the role of the optimisation is to identify a spatial intervention that improves it. In the setting considered here, this can occur by modifying either the bifurcation structure of the system or the realised trajectory followed by the state. \Cref{fig:FIGURE1} provides an overview of this idea: it shows a system represented by \eqref{eq:AC2} (left) and its optimised version (right). It illustrates how a targeted intervention can shift a critical transition to more adverse environmental conditions, thereby avoiding a tipping event over the parameter range considered. The figure is intended as a conceptual guide for the optimisation problems introduced below.
In this work, we distinguish two optimisation settings: equilibrium optimisation, which aims to modify the bifurcation structure to delay tipping, and time-dependent (final-state) optimisation, which influences transient spatial dynamics, including local tipping and front propagation.
 In equilibrium optimisation, resilience is assessed from stationary solutions. The intervention $z(x)$ is fixed, and the state satisfies the steady-state equations, imposed as constraints.
 By contrast, in time-dependent optimisation, resilience is quantified from system trajectories. The model is integrated forward in time from a prescribed initial condition over a chosen time horizon, under a fixed spatial intervention $z(x)$. The governing dynamics enter as equality constraints, ensuring that each candidate intervention generates a valid trajectory. The objective function evaluates resilience along the trajectory or at its final state. 
This formulation is closely related to the conditional nonlinear optimal perturbation (CNOP) framework, where perturbations are optimised to maximise an objective after finite-time integration \cite{Mu2003_CNOP,Oosterwijk2017_CNOP}, a common approach in nonlinear stability and predictability studies. In our case, the perturbation corresponds to a spatial intervention on the system parameters rather than an initial-condition perturbation.

In both optimisation settings, the intervention $z(x)$ is defined on a continuous spatial domain, so the admissible set $\Omega$ is, in principle, infinite-dimensional. 
In practice, optimisation is carried out over a finite-dimensional parametrisation of $z(x)$. This discretisation  is not unique: depending on the structure of the system under study, one could employ, for example, Chebyshev polynomials, Fourier modes, or Lagrange polynomials. The particular parametrisation adopted in this paper is introduced together with the Allen--Cahn application in \Cref{sec:4}.

In the following section, we apply this framework to a spatially extended gradient system described by the Allen--Cahn equation.

\section{Application of the Optimisation Framework to the Allen-Cahn Equation}\label{sec:4}

We apply the resilience optimisation framework to the one-dimensional Allen--Cahn equation \cite{bastiaansen2022_fragmented_tipping, allen1972_allen_cahn, Bray1994_AllenCahn} as a minimal conceptual model representative of bistable gradient systems:
\begin{equation}\label{eq:AC_general}
    \frac{\partial y}{\partial t}
    = \frac{\partial^2 y}{\partial x^2}
    + y\left(1 - y^2\right)
    + \mu(t) + \mu_{\mathrm{het}}(x)
    + \mu_{\mathrm{pert}}(x)\, .
\end{equation}
Here, $\mu(t)$ is a scalar bifurcation parameter, $\mu_{\mathrm{het}}(x)$ is a pre-existing spatial heterogeneity fixed prior to any intervention, and $\mu_{\mathrm{pert}}(x)$ is an added spatial perturbation that acts as the decision variable in the optimisation, i.e. $z(x) = \mu_{\mathrm{pert}}(x)$.

In the spatially uniform setting, where $\mu_{\mathrm{het}}(x) = \mu_{\mathrm{pert}}(x) = 0$, the system exhibits two saddle-node bifurcations at
\[
\mu_{\mathrm{tip}}^{\pm} = \pm \frac{2}{9}\sqrt{3}.
\]
For $\mu_{\mathrm{tip}}^- < \mu < \mu_{\mathrm{tip}}^+$, the system is bistable, with two competing homogeneous equilibria $y^{\pm}$. The Maxwell point is located at $\mu_M = 0$ and marks the parameter value at which an isolated front connecting the two states has zero speed (see \Cref{sec:2}). For $\mu \neq 0$, the sign of $\mu$ biases the local potential and determines the direction of front propagation.
When spatial heterogeneity is present, i.e. $\mu_{\mathrm{het}}(x)$ is non-constant, translation invariance is broken, and the equilibrium structure becomes richer \cite{bastiaansen2022_fragmented_tipping} (see \Cref{fig:0}). Steady states may include pinned fronts, and numerical continuation in $\mu$ can reveal multiple saddle-node bifurcations creating new equilibrium branches. We denote these saddle-nodes by $\mu_{\mathrm{SN},k}$, where the index $k$ labels the successive folds induced by the spatial heterogeneity.

In the following, increasing values of the bifurcation parameter $\mu$ are interpreted as increasing external pressure on the system. We label the lower homogeneous state as desirable and the upper homogeneous state unfavourable. This convention fixes the interpretation of tipping, recovery, and resilience in all optimisation problems considered below.
Resilience is quantified using problem-specific proxies tailored to the optimisation setting. 
In the equilibrium optimisation, resilience is associated with the saddle-node value $\mu_{\mathrm{SN},k}$ at which lower equilibria (locally) cease to exist, and the optimisation seeks perturbations that shift a selected saddle-node to higher values of $\mu$. In the final-state optimisations, resilience is measured from system trajectories, using the spatial mean of the solution at a prescribed final time as a proxy for the extent of the desirable state. The optimisation seeks perturbations that minimise this quantity. 

Within this setting, we study three optimisation problems for the Allen--Cahn equation: one equilibrium optimisation and two final-state optimisations. In the latter, the bifurcation parameter is either fixed in time (off-branch optimisation) or slowly varying so that trajectories track stable branches (on-branch optimisation).

\begin{figure}
    \centering
    \includegraphics[width=\linewidth]{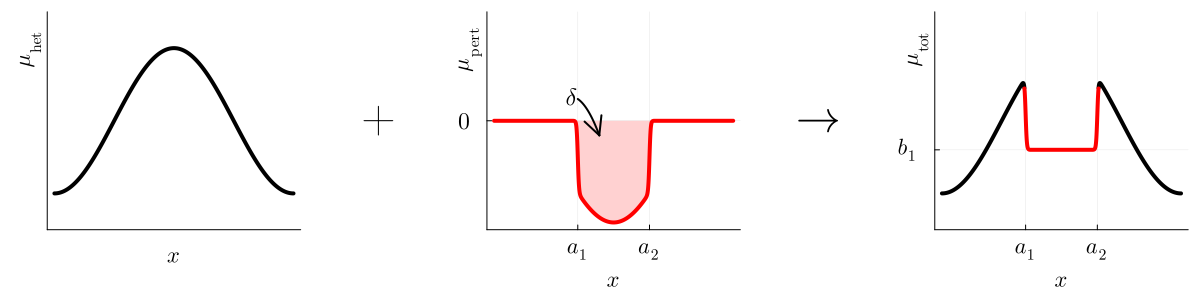}
    \caption{Schematic illustration of the perturbation parametrisation for a single local modification ($N=1$). \textbf{Left:} pre-existing heterogeneity profile $\mu_{\mathrm{het}}(x)$ (black). \textbf{Centre:} perturbation $\mu_{\mathrm{pert}}(x)$ (red), constructed from one interval function acting on the interval $(a_1,a_2)$. The shaded red area is proportional to the perturbation size $\delta$, corresponding to the size constraint imposed on the intervention in the optimisation problems below. \textbf{Right:} resulting total heterogeneity $\mu_{\mathrm{tot}}(x)=\mu_{\mathrm{het}}(x)+\mu_{\mathrm{pert}}(x)$. Outside the interval $(a_1,a_2)$, the total heterogeneity coincides with the pre-existing profile (black), while inside the interval it is locally flattened towards the target level $b_1$ (red). The interval endpoints $a_1$ and $a_2$ are marked on the horizontal axis, and the target plateau level $b$ is marked on the vertical axis.}
    \label{fig:parametrisation}
\end{figure}

To make the optimisation problem finite-dimensional, we parametrise the intervention $\mu_{\mathrm{pert}}(x)$ using a finite family of local modifications of the pre-existing heterogeneity $\mu_{\mathrm{het}}(x)$. The idea is to let the optimisation act by selecting spatial intervals and locally flattening $\mu_{\mathrm{het}}(x)$ onto prescribed constant levels. In this way, the intervention directly controls where the total heterogeneity is pushed relative to key thresholds such as the saddle-node values and the Maxwell point, while keeping the search space low-dimensional and interpretable. This construction is illustrated schematically in \Cref{fig:parametrisation}.
To implement this idea, we introduce a basic interval function
\[
g(x;a_1,a_2)=
\frac{1}{4}
\Bigl(\tanh\bigl(\kappa(x-a_1)\bigr)+1\Bigr)
\Bigl(\tanh\bigl(\kappa(x-a_2)\bigr)-1\Bigr),
\]
with $a_1<a_2$. The parameter $\kappa$ is set equal to $100$, so $g(x;a_1,a_2)$ is approximately equal to $-1$ on the interval $(a_1,a_2)$ and approximately $0$ outside of it. Using this function, we define a local perturbation term
\[
\bar g(x;a_1,a_2,b)=g(x;a_1,a_2)\bigl(\mu_\mathrm{het}(x)-b\bigr),
\]
where $b$ is a prescribed target level. Since $g(x;a_1,a_2)\approx -1$ on $(a_1,a_2)$, adding $\bar g$ to $\mu_{\mathrm{het}}(x)$, locally replaces the original heterogeneity profile by an approximately constant plateau $b$, while leaving it essentially unchanged outside the chosen interval.
The full perturbation is then constructed as a sum of such local contributions,
\[
\mu_{\mathrm{pert}}(x)=\sum_{k=1}^{N} \bar g(x;a_{2k-1},a_{2k},b_k),
\]
with $a_{2k-1} < a_{2k}$, so that the total heterogeneity becomes
\[
\mu_{\mathrm{tot}}(x):=\mu + \mu_{\mathrm{het}}(x)+\mu_{\mathrm{pert}}(x).
\]
The optimisation is therefore performed over the finite set of parameters $\{(a_{2k-1}, a_{2k}, b_k)\}_{k=1}^N$ which determine the placement and strength of the local flattening. As shown in \Cref{fig:parametrisation}, this parametrisation allows the optimisation to act directly on the local mechanisms governing tipping and front pinning, while remaining transparent and computationally manageable.

All optimisation problems are solved using a gradient-based sequential quadratic programming method (SLSQP) \cite{SLSQP}, as implemented in NLopt \cite{NLopt}, with gradients computed via automatic differentiation in Julia \cite{Julia2017}. Bifurcation diagrams are computed using numerical continuation methods as implemented in BifurcationKit.jl \cite{BifKit}. The code used to produce these results is available at the following link \url{https://github.com/AuroraFaureRagani/optimisation_spatially_extended_systems}.

\subsection{Equilibrium optimisation}
The first application of the framework concerns equilibrium optimisation, aimed at delaying tipping by modifying the steady states of the Allen–Cahn equation. To study this setting, we consider the stationary form of \eqref{eq:AC_general}, obtained by setting the time derivative to zero and fixing the bifurcation parameter in time ($\mu(t)=\mu$): 
\begin{equation}
    0 =\frac{\partial^2 y}{\partial x^2} + y\left(1 - y^2\right) + \mu + \mu_{\mathrm{het}}(x) + \mu_{\mathrm{pert}}(x)
\end{equation}
This equation acts as the state constraint in the optimisation problem. As introduced above, we denote the $k^\text{th}$ saddle-node by $\mu_{\mathrm{SN},k}$.
In this setting, resilience is quantified by this parameter value $\mu_{\mathrm{SN},k}$ at which the lower-state equilibrium ceases to exist locally through a saddle-node bifurcation. The optimisation therefore seeks a spatial perturbation $\mu_{\mathrm{pert}}(x)$ that shifts a chosen tipping point to the highest possible value of $\mu$, thus postponing tipping as much as possible. Accordingly, the objective function is taken as the negative of the tipping value $\mu_{\mathrm{SN},k}$. The constraints enforce that the solution $y$ is a steady state located at a saddle-node and that the perturbation satisfies a size constraint, expressed as an $L^1$-bound that limits the total magnitude of the intervention. In \Cref{fig:parametrisation}, this is represented schematically by the shaded area in the centre panel.
Consequently, we seek $\mu_{\mathrm{pert}}(x)$ that solves the following:
\begin{align*}
    & \min_{\mu_{\mathrm{pert}}}  (-\mu_{\mathrm{SN},k})  \\
    &\text{subject to } 
    \begin{cases}
        0 = \dfrac{\partial^2 y}{\partial x^2} + y(1 - y^2) + \mu_{\mathrm{SN},k} + \mu_{\mathrm{het}}(x) + \mu_{\mathrm{pert}}(x), \\[0.4em]
        \text{$(y,\mu_{\mathrm{SN},k})$ is a saddle-node steady state},\\[0.4em]
        \displaystyle \frac{1}{2}\int_{-1}^1 |\mu_{\mathrm{pert}}(x)|\,dx \le \delta.
    \end{cases}
\end{align*}
The saddle-node condition is imposed by requiring that $y$ is a steady state and that the linearisation of the steady-state equation at $(y,\mu_{\mathrm{SN},k})$ has a vanishing dominant eigenvalue. In the numerical implementation, this condition is enforced by requiring that the largest eigenvalue of the Jacobian of the discretised system vanishes. 

To illustrate the equilibrium optimisation, we consider two cases. 
Case A uses the pre-existing heterogeneity $\mu_{\mathrm{het}}(x) = \tfrac{1}{2}\cos(\pi x)$ and targets the first saddle node ($k=1$), 
while Case A considers $\mu_{\mathrm{het}}(x) = \tfrac{1}{2}\cos(\pi x)\sin(2\pi x)$ and targets the fifth saddle node ($k=5$). We interpret $\delta$ as a budget on the intervention size. We fix it at $\delta = 0.3$ for both cases so that  $\mu_{\mathrm{pert}}(x)$ can meaningfully reshape  $\mu_{\mathrm{het}}(x)$ while remaining comparable to its amplitude.
The corresponding results are shown in \Cref{fig:1A} for Case A and \Cref{fig:1B} for Case B.

\subsubsection{Case A: Delayed full tipping}
\begin{figure}
    \centering
    \includegraphics[width=\linewidth]{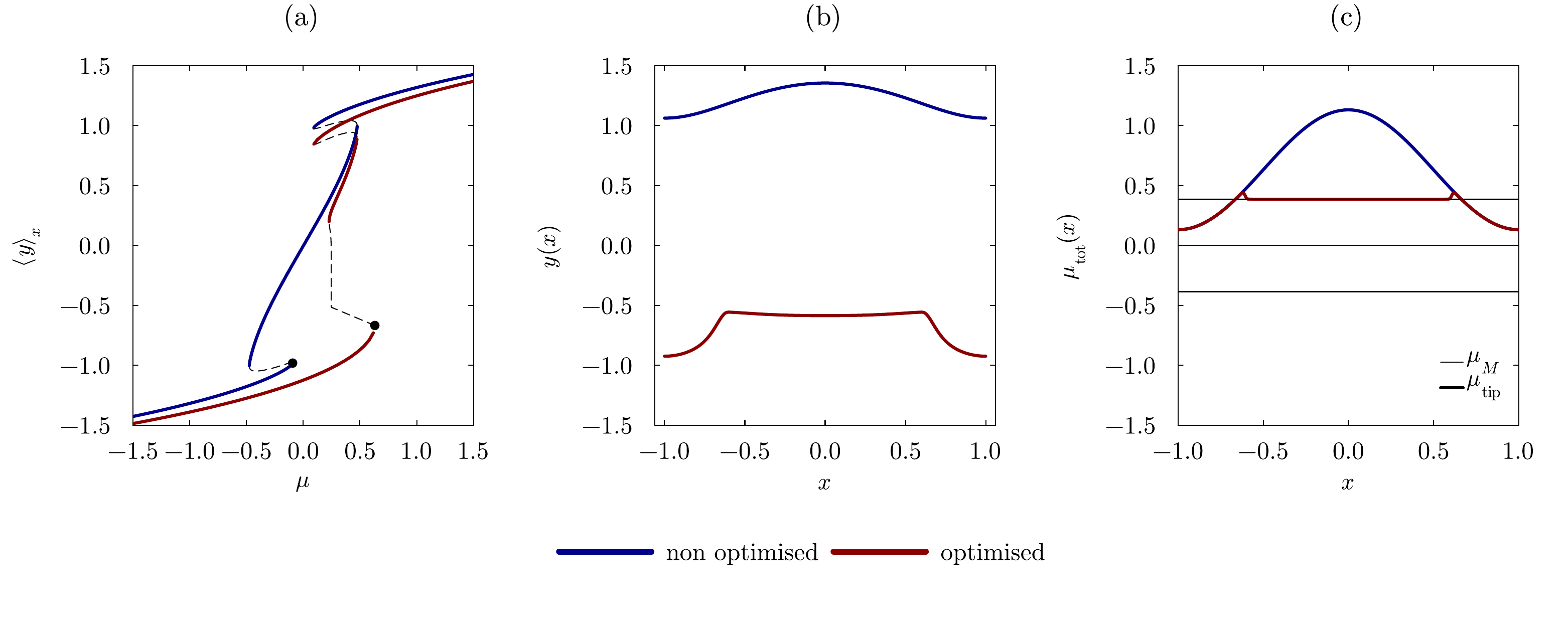}
    \caption{
Equilibrium optimisation targeting the location of the first saddle-node $\mu^{\mathrm{opt}}_{\mathrm{SN},1}$, with heterogeneity $\mu_{\mathrm{het}}(x) = \tfrac{1}{2}\cos(\pi x)$.
\textbf{(a)} Bifurcation diagrams showing the spatial mean of the solution, $\langle y  \rangle_x = \frac{1}{2}\int_{-1}^{1}y(x)dx$, as a function of the bifurcation parameter $\mu$. Thick curves denote stable branches, while dashed curves denote unstable branches. Plots are shown for the non-optimised system (dark blue) and the optimised system (dark red). Black dots indicate the targeted saddle-node in both diagrams.
\textbf{(b)} Spatial profiles $y(x)$ corresponding to the optimal parameter value $\mu = \mu^{\mathrm{opt}}_{\mathrm{SN},1}$ for the non-optimised (dark blue) and optimised (dark red) systems.
\textbf{(c)} Spatial heterogeneities $\mu_{\mathrm{tot}}(x)$ with $\mu = \mu^{\mathrm{opt}}_{\mathrm{SN},1}$ for the non-optimised system (dark blue) and the optimised system (dark red). Horizontal lines indicate the local tipping points $\mu_{\mathrm{tip}}^{\pm}$ (thick line) and the Maxwell point $\mu_M$ (thin line).
}

    \label{fig:1A}
\end{figure}

We begin by analysing Case A, corresponding to the heterogeneity $\mu_{\mathrm{het}}(x) = \tfrac{1}{2}\cos(\pi x)$, to assess how the optimisation shifts the first tipping point. \Cref{fig:1A} illustrates the effects of this optimisation.
Panel (a) shows the bifurcation diagrams of the non-optimised (blue) and optimised (red) systems. In the non-optimised case, the tipping point occurs at $\mu_{\mathrm{SN},1} = -0.09$, whereas optimisation shifts it to $\mu_{\mathrm{SN},1}^{\mathrm{opt}} = 0.63$. This latter value lies beyond the multistable regime of the non-optimised system. Outside that regime, the only possible stable solution is the fully tipped state. Therefore, the optimisation significantly enhances resilience by delaying tipping to larger values of $\mu$.

\Cref{fig:1A}(c) shows the total heterogeneity $\mu_{\mathrm{tot}}(x)$ defined as
\[
\mu_{\mathrm{tot}}(x; \mu)= \mu + \mu_{\mathrm{het}}(x) + \mu_{\mathrm{pert}}(x),
\]
where the bifurcation parameter is at the saddle node of the optimal case ($\mu = \mu_{\mathrm{SN},1}^{\mathrm{opt}}$).
In this plot, the blue and red lines represent the non-optimised ($\mu_{\mathrm{pert}}(x)$ = 0) and optimised case, respectively. Although $\mu_{\mathrm{pert}}(x)$ is not plotted directly, its effect can be inferred from the difference between the two total heterogeneity profiles. The optimisation modifies the pre-existing heterogeneity so that $\mu_{\mathrm{tot}}(x)$ remains below $\mu_{\mathrm{tip}}^+$, except on a narrow interval. For the diffusion levels used here, this interval is insufficient to trigger local tipping. The optimal perturbation $\mu_{\mathrm{pert}}(x)$ is non-zero only where the non-optimised $\mu_{\mathrm{tot}}(x)$ exceeds $\mu_{\mathrm{tip}}^+$. Within this region ($x \in (-0.61, 0.61)$), it acts to keep the local heterogeneity slightly below $\mu_{\mathrm{tip}}^+$ to prevent tipping. Outside this interval, $\mu_{\mathrm{pert}}(x)$ vanishes, and the total heterogeneity coincides with the non-optimised profile. 

At $\mu = \mu_{\mathrm{SN},1}^{\mathrm{opt}}$, the outer regions of the domain are in the undesirable state in the non-optimised case, whereas they remain untipped in the optimised case (\Cref{fig:1A}(b)), despite identical local heterogeneity there.
This distinction arises from the presence of propagating fronts in the non-optimised system. Once the centre of the domain tips, two fronts form and move outward. Their propagation continues until they reach points where $\mu_{\mathrm{tot}}(x)$ equals the Maxwell point, at which the fronts become pinned and can no longer advance. As $\mu$ increases, these fronts in the non-optimised case reach the boundaries of the domain and disappear, getting the solution to a fully tipped state at parameter values lower than $\mu_{\mathrm{SN},1}^{\mathrm{opt}}$. In the optimised case, no front is generated; keeping $\mu_{\mathrm{tot}}(x)$ below $\mu_{\mathrm{tip}}^+$ everywhere is sufficient to prevent tipping.

\subsubsection{Case B: Delayed fragmented tipping}
\begin{figure}
    \centering
    \includegraphics[width=\linewidth]{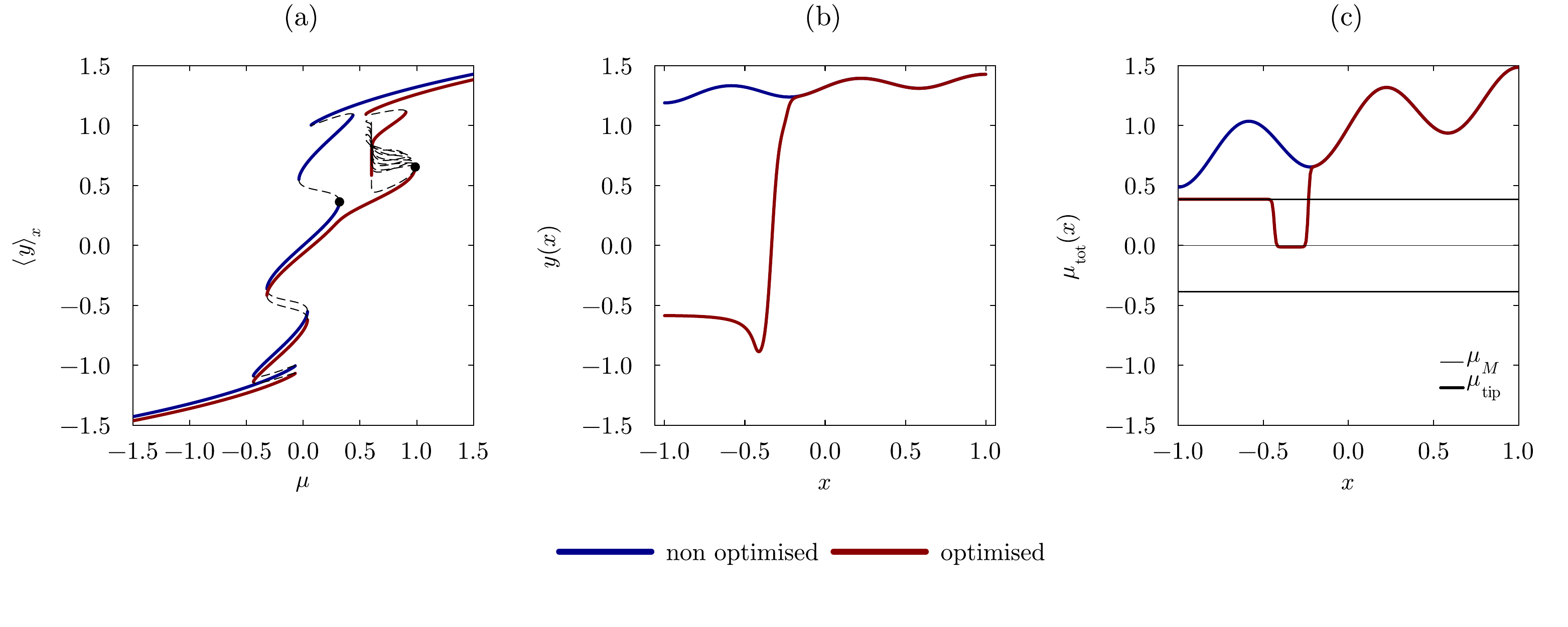}
    \caption{
Equilibrium optimisation targeting the location of the fifth saddle-node $\mu^{\mathrm{opt}}_{\mathrm{SN},5}$, with heterogeneity $\mu_{\mathrm{het}}(x) = \tfrac{1}{2}\cos(\pi x) \sin(2 \pi x)$.
\textbf{(a)} Bifurcation diagrams showing the spatial mean of the solution, $\langle y  \rangle_x = \frac{1}{2}\int_{-1}^{1}y(x)dx$, as a function of the bifurcation parameter $\mu$. Thick curves denote stable branches, while dashed curves denote unstable branches. Plots are shown for the non-optimised system (dark blue) and the optimised system (dark red). Black dots indicate the targeted saddle-node in both diagrams.
\textbf{(b)} Spatial profiles $y(x)$ corresponding to the optimal parameter value $\mu = \mu^{\mathrm{opt}}_{\mathrm{SN},5}$ for the non-optimised (dark blue) and optimised (dark red) systems.
\textbf{(c)} Spatial heterogeneities $\mu_{\mathrm{tot}}(x)$ with $\mu = \mu^{\mathrm{opt}}_{\mathrm{SN},5}$ for the non-optimised system (dark blue) and the optimised system (dark red). Horizontal lines indicate the local tipping points $\mu_{\mathrm{tip}}^{\pm}$ (thick line) and the Maxwell point $\mu_M$ (thin line).
}
    \label{fig:1B}
\end{figure}
To test the framework on a more complex heterogeneity, we consider the equilibrium optimisation of Case B, where the pre-existing heterogeneity is $\mu_{\mathrm{het}}(x) = \tfrac{1}{2} \cos(\pi x) \sin(2 \pi x)$. In this case, the non-optimised problem has eight saddle nodes, and the optimisation aims to postpone the fifth one, counting from below. 
\Cref{fig:1B}(a) shows the non-optimised (in blue) and optimised (in red) bifurcation diagrams. The optimisation shifts the fifth saddle node from $\mu_{\mathrm{SN},5} = 0.32$ to $\mu_{\mathrm{SN},5}^{\mathrm{opt}} = 1.0$. In panel (c), the blue and red lines represent the non-optimised and optimised total heterogeneities $\mu_{\mathrm{tot}}(x)$, respectively, both evaluated at $\mu = \mu_{\mathrm{SN},5}^{\mathrm{opt}}$. The non-optimised heterogeneity lies entirely above the homogeneous tipping point, corresponding to a fully tipped steady state where no front can persist. In contrast, the optimisation alters the total heterogeneity so that the system remains in a partially tipped steady state at $\mu_{\mathrm{SN},5}^{\mathrm{opt}}$. 

At the fifth saddle node, the solution is tipped in the right portion of the domain, while the left part remains untipped. This partially tipped situation persists when $\mu$ increases, provided two conditions are satisfied. First, in the regions that should remain untipped, the total heterogeneity must remain below the homogeneous tipping point; otherwise, local tipping would occur. Second, the existing front must be pinned to prevent further propagation. The optimisation achieves both requirements. It acts only on the left region of the domain where the solution has not tipped at the non-optimised fifth saddle node. It introduces a narrow interval in which $\mu_{\mathrm{tot}}(x) < \mu_M$ so that $\mu_{\mathrm{tot}}(x)$ crosses the Maxwell point at the boundaries of this interval. To the left of this narrow region, the total heterogeneity is modified to be slightly below the homogeneous tipping point $\mu_{\mathrm{tip}}^+$.

\subsection{Off-branch final-state optimisation}
The second type of optimisation applied to the Allen-Cahn equation is final-state optimisation. Here, we consider solutions that evolve in time from a prescribed initial condition at a fixed value of the bifurcation parameter $\mu$. 
We fix $\mu$ within the bistable regime. Tipping and front dynamics are then determined by the total heterogeneity $\mu_{\mathrm{tot}}(x) = \mu + \mu_{\mathrm{het}}(x)$, and can be modified through the applied perturbation $\mu_{\mathrm{pert}}(x)$. Because the time evolution is not necessarily close to a steady state, this optimisation is referred to as off-branch optimisation. The goal in this optimisation problem is to find the optimal spatial perturbation that minimises the mean value of the solution at a prescribed final time $t_e$, which serves as a scalar proxy for the spatial extent of the unfavourable state. 
Here, $\langle \cdot\rangle_x$ denotes the spatial average over the domain, i.e.
\[
\langle y(x,t_e)\rangle_x = \frac{1}{2}\int_{-1}^{1} y(x,t_e)\,dx.
\]
In this setting, the optimal perturbation can influence both the creation or disruption of fronts and the direction of front migration mechanisms that ultimately determine the system final state. Formally, we seek $\mu_{\mathrm{pert}}$ such that
\begin{align*}
    &\min_{\mu_{\mathrm{pert}}} \langle y(x, t_e)\rangle_x \\
    &\text{subject to } 
    \begin{cases}
        \dfrac{\partial y}{\partial t}= \dfrac{\partial^2 y}{\partial x^2} + y(1 - y^2) + \mu + \mu_{\mathrm{het}}(x) + \mu_{\mathrm{pert}}(x),  \qquad t \in [t_0, t_e]\\[0.4em]
        y(x, t_0) = y_{\mathrm{IC}},\\[0.4em]
        \displaystyle \frac{1}{2}\int_{-1}^1 |\mu_{\mathrm{pert}}(x)|\,dx \le \delta.
    \end{cases}
\end{align*}
Here, the first two constraints ensure the Allen-Cahn dynamics with given initial condition $y_{\mathrm{IC}}$, while the last one imposes the size constraint on the perturbation. The budget $\delta$ limits the size of  $\mu_{\mathrm{pert}}(x)$. We set $t_e$ so that fronts have either pinned or reached the boundary, so $\langle y(x, t_e)\rangle_x$ displays the equilibrium rather than the transient.
In the simulations that we analyse, two different initial conditions are considered. The first one (case A) contains an up-front at the domain centre, $y_{\mathrm{IC}}^A(x) = \tanh(100x)$, and optimisation halts the front and promotes recovery through favourable front propagation, leading to expansion of the lower state. In Case B, the system starts from a fully tipped state, $y_{\mathrm{IC}}^B(x) = 1$, and optimisation instead induces local front creation. Recovery in this case occurs only if the resulting fronts subsequently propagate so as to expand the lower state. For both cases, we test two values of $\delta$, with and without a pre-existing heterogeneity.

\subsubsection{Case A: Front pinning and state recovery}
\begin{figure}
    \centering
    \includegraphics[width=1\linewidth]{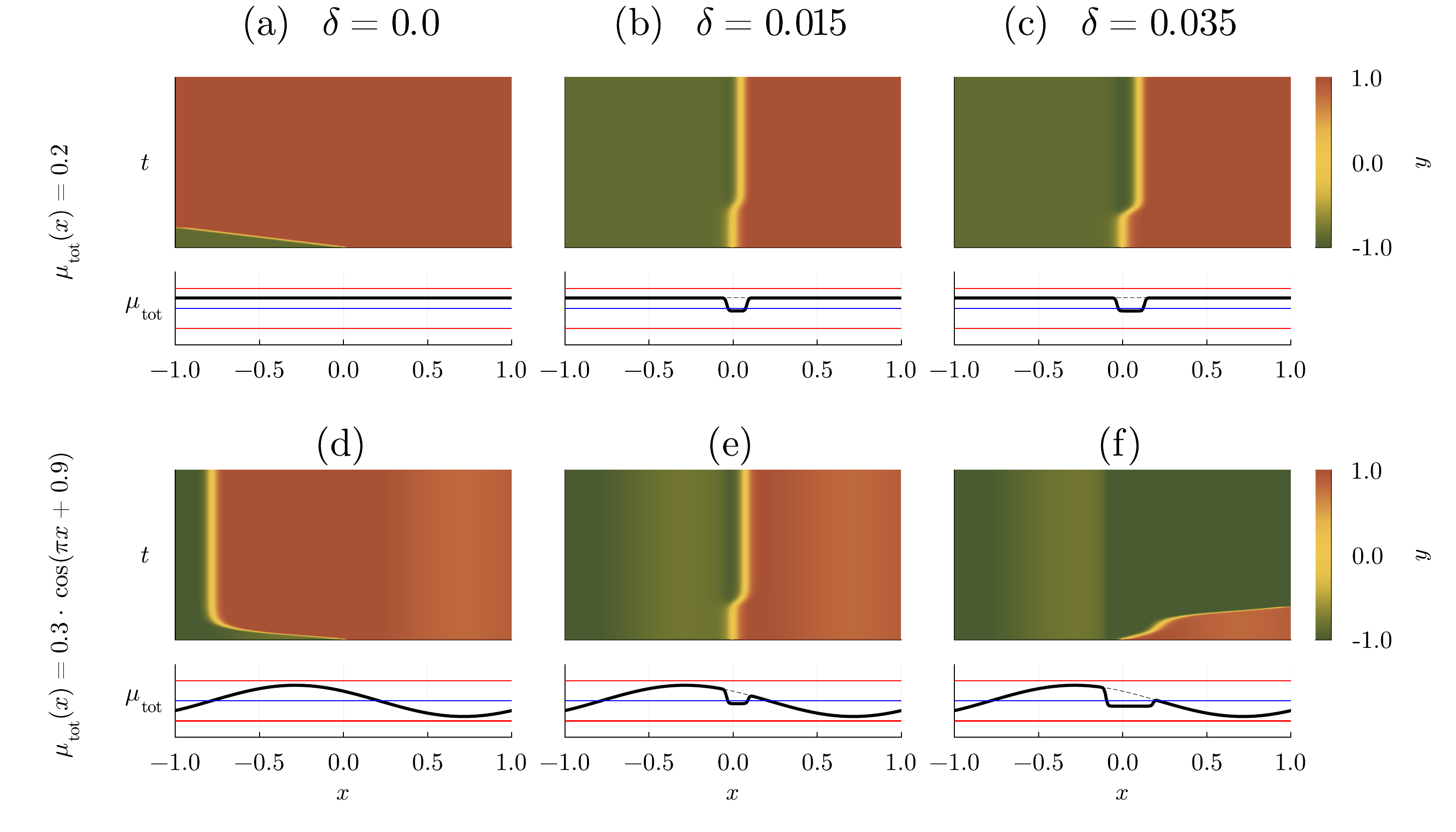}
        \caption{
Off-branch final-state optimisation starting from a front-like initial
condition $y_{\mathrm{IC}}(x)=\tanh(100x)$ (Case A).
\textbf{(a)--(c)} Homogeneous case with $\mu_{\mathrm{het}}(x)=0$ and
$\mu=0.2$.
\textbf{(d)--(f)} Heterogeneous case with $\mu=0$ and
$\mu_{\mathrm{het}}(x)=0.3\cos(\pi x+0.9)$.
From left to right, the perturbation size is $\delta=0$ in panels
(a) and (d), $\delta=0.015$ in panels (b) and (e), and
$\delta=0.035$ in panels (c) and (f).
Each panel consists of two subplots: the upper subplot shows the
space--time evolution of $y(x,t)$, and the lower subplot shows the
associated total heterogeneity $\mu_{\mathrm{tot}}(x)$, with horizontal
lines indicating the local tipping thresholds $\mu_{\mathrm{tip}}^\pm$
(red) and the Maxwell point $\mu_M$ (blue).
}\label{fig:2A}
\end{figure}

We first examine the case where the initial condition has a front located at $x=0$. The results for different heterogeneities $\mu_{\mathrm{het}}(x)$ and different perturbations  $\mu_{\mathrm{pert}}(x)$ are shown in \Cref{fig:2A}. Each panel has the time evolution of the solution above and the associated total heterogeneity $\mu_{\mathrm{tot}}(x) = \mu + \mu_{\mathrm{het}}(x) + \mu_{\mathrm{pert}}(x)$ below.
Panels (a)--(c) show the homogeneous case ($\mu_{\mathrm{het}}(x) = 0$) with the bifurcation parameter fixed at $\mu = 0.2$, so that front motion, in the absence of any applied perturbation, is driven by a uniform bias.
Panels (d)--(f) show the heterogeneous case, with $\mu=0$ and a sinusoidal heterogeneity ($\mu_{\mathrm{het}}(x) = 0.3 \cos(\pi x + 0.9)$), so that front motion is determined by spatial variations in $\mu_{\mathrm{tot}}(x)$.
From left to right, the perturbation size increases across the columns: panels (a) and (d) show the base case without optimisation ($\delta=0$), panels (b) and (e) show the optimised case with $\delta=0.015$, and panels (c) and (f) show the optimised case with $\delta=0.035$. These two $\delta$ values have been chosen to illustrate a low-budget regime that triggers only front halting and a higher-budget regime that also sustains favourable front propagation and larger recovery.

We first analyse how the front dynamics evolve in Case A under increasing perturbation sizes and in the presence or absence of heterogeneity.
In the homogeneous case, shown in \Cref{fig:2A}(a)--(c), no pre-existing heterogeneity is present. In the base case, without optimisation ($\delta = 0$, panel (a)), the initial up-front propagates to the right because $\mu_{\mathrm{tot}}(x) = \mu = 0.2$ exceeds the Maxwell point everywhere, leading the system to the fully tipped state,  meaning that the upper state occupies the entire domain, with spatial mean $\langle y(x, t_e) \rangle_x = 1$. Allowing a small-sized perturbation ($\delta = 0.015$, panel (b)) substantially alters the outcome: the mean of the final state nearly halves ($\langle y(x, t_e) \rangle_x = 0.05$). This occurs because the optimal perturbation lowers $\mu_{\mathrm{tot}}(x)$ below the Maxwell point over a very narrow interval around the front location, stopping the front propagation and in fact inducing a short-lived motion of the front to the right. However, $\mu_{\mathrm{tot}}(x)$ exceeds $\mu_M$ again immediately outside this region, bringing the front velocity to zero. As a result, the front experiences only a slight displacement before becoming effectively pinned. For a larger perturbation size ($\delta = 0.035$, panel (c)), the region where $\mu_{\mathrm{tot}}(x)$ falls below the Maxwell point can widen, and the front begins to retreat. Consequently, the mean of the final state decreases further ($\langle y(x, t_e) \rangle_x = -0.11$), indicating partial recovery of the untipped state compared to the initial condition.

We now turn to the heterogeneous case, shown in \Cref{fig:2A}(d)--(f), where a pre-existing heterogeneity is $\mu_{\mathrm{het}}(x) = 0.3 \cos(\pi x + 0.9)$. This choice provides a smooth spatial variation that remains within the bistable regime while introducing regions above and below the Maxwell point, thereby enabling both favourable and unfavourable front motion.
In the base case without optimisation ($\delta = 0$, panel (d)), the front moves leftward wherever $\mu_{\mathrm{tot}}(x)$ lies below the Maxwell point, and the system approaches a partially tipped state with spatial mean $\langle y(x, t_e) \rangle_x = 0.76$. Applying optimisation, however, can counteract this shift. With a small perturbation ($\delta = 0.015$, panel (e)), the front can be halted, and the mean of the final state decreases to $\langle y(x, t_e) \rangle_x = -0.08$. For a larger perturbation size ($\delta = 0.035$, panel (f)), the optimal perturbation enables complete recovery of the untipped state ($\langle y(x, t_e) \rangle_x = -1.0$). Here, the heterogeneity helps with the state recovery, as the optimisation only needs to lower $\mu_{\mathrm{tot}}(x)$ slightly below the Maxwell point between the front position and the zero of the heterogeneity. As a result, the front retreats to the right and the system returns to an untipped state. In this specific setting, heterogeneity makes the optimisation more effective, allowing full recovery with the same perturbation size that only achieves partial recovery in the homogeneous case.

\subsubsection{Case B: Front creation and state recovery}
\begin{figure}
    \centering
    \includegraphics[width=\linewidth]{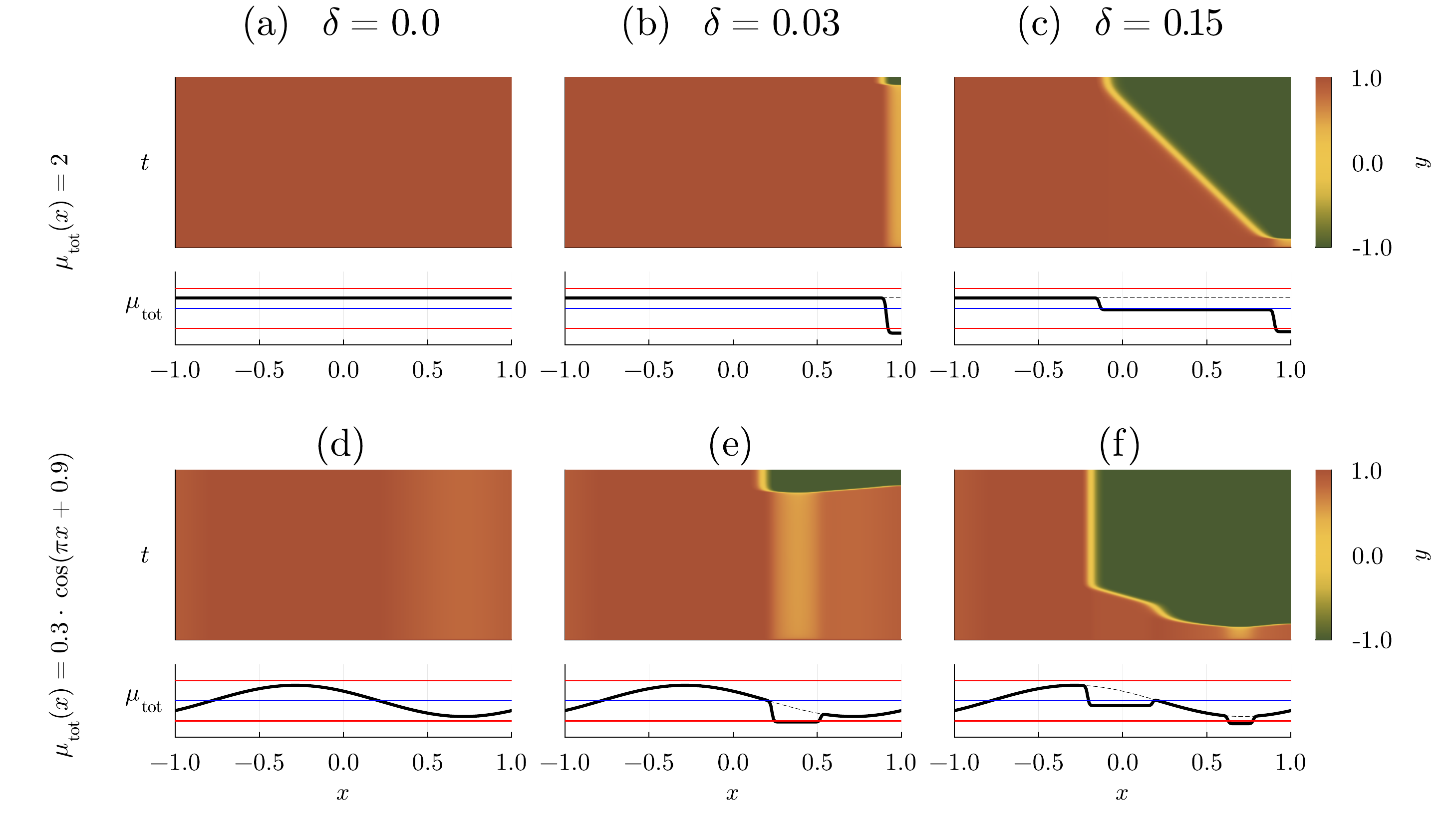}
        \caption{
Off-branch final-state optimisation starting from a fully tipped initial
condition $y_{\mathrm{IC}}(x)=1$ (Case B).
\textbf{(a)--(c)} Homogeneous case with $\mu_{\mathrm{het}}(x)=0$ and
$\mu=0.2$.
\textbf{(d)--(f)} Heterogeneous case with $\mu=0$ and
$\mu_{\mathrm{het}}(x)=0.3\cos(\pi x+0.9)$.
From left to right, the perturbation size is $\delta=0$ in panels
(a) and (d), $\delta=0.03$ in panels (b) and (e), and
$\delta=0.15$ in panels (c) and (f).
Each panel consists of two subplots: the upper subplot shows the space--time evolution of $y(x,t)$, and the lower subplot shows the associated total heterogeneity $\mu_{\mathrm{tot}}(x)$, with horizontal lines indicating the local tipping thresholds $\mu_{\mathrm{tip}}^\pm$
(red) and the Maxwell point $\mu_M$ (blue).
}

    \label{fig:2B}
\end{figure}

Unlike Case A, where we halted or reversed a pre-existing front, Case B begins from a fully tipped state and therefore must first create front(s) before their propagation can be controlled.
The results for Case B are shown in \Cref{fig:2B}. Panels (a)--(c) show the homogeneous case ($\mu_{\mathrm{het}}(x)=0$), while panels (d)--(f) show the heterogeneous case with $\mu_{\mathrm{het}}(x)=0.3\cos(\pi x+0.9)$.  From left to right, the perturbation size is $\delta=0$, $\delta=0.03$, and $\delta=0.15$.

In both base cases (with and without heterogeneity), the state remains tipped throughout the simulation if no optimisation is performed (panels (a) and (d)). This persistence of the tipped state occurs because the total heterogeneity stays above the lower homogeneous tipping threshold, preventing any part of the domain from recovering.
The solution behaviour changes once optimisation is introduced (panels (b) and (e)). Even a relatively small perturbation can substantially modify the system time evolution. The optimal perturbation works by locally reducing the total heterogeneity, $\mu_{\mathrm{tot}}(x)$, below the homogeneous tipping point $\mu_{\mathrm{tip}}^-$, so that the lower state can nucleate; subsequent front propagation is then controlled by the Maxwell point. The location and extent of this recovery depend on the shape of the perturbation $\mu_{\mathrm{pert}}(x)$, and therefore differ between the homogeneous and heterogeneous settings. For instance, in \Cref{fig:2B}(b), the perturbation acts near the domain boundary, generating a single front. In contrast, when $\mu_{\mathrm{het}} = 0.3 \cos (  \pi x + 0.9 )$ in \Cref{fig:2B}(e), the optimisation targets an internal region where $\mu_{\mathrm{het}}(x)<\mu_M$. Because this region is not located at the boundary, the solution recovers in the middle of the domain, creating two fronts that travel in opposite directions. These fronts move outward until they reach the points where the total heterogeneity becomes zero. In this case, the heterogeneity facilitates recovery: although the perturbation acts only on a small part of the domain, it triggers front propagation of roughly half of it ($\langle y(x, t_e) \rangle_x = 0.16$). Because the local heterogeneity lies below the Maxwell point, the resulting front becomes unpinned and can propagate freely, allowing the desirable state to spread into regions where no perturbation is applied.
 
When the perturbation size is increased (panels (c) and (f)), the optimisation not only creates fronts but also enhances recovery by unpinning them wherever possible. In the initially homogeneous case (panel (c)), the front first forms where the total heterogeneity is locally below the homogeneous tipping point. Thereafter, it is sufficient that the total heterogeneity remains below the Maxwell point for it to propagate. The extent of the region where $\mu_{\mathrm{tot}}(x) < \mu_M$ is determined by the perturbation size constraint; once the front reaches a point where this condition no longer holds, it becomes pinned and stops moving. In this case, the final state has mean $\langle y(x, t_e) \rangle_x = -0.07$.
 In the heterogeneous case (panel (f)), a single interior recovery region is again formed, producing two fronts that move outward without requiring further intervention. The right-moving front reaches the domain boundary and disappears. The left-moving front, however, would stall at the location where $\mu_{\mathrm{tot}}(x) = \mu_M$ if no additional action were taken. Because the perturbation budget is larger than in the previous case, the optimisation can also lower the pre-existing heterogeneity below the Maxwell point in that region, allowing the left-moving front to continue travelling and complete the recovery of the domain ($\langle y(x, t_e) \rangle_x = -1.06$). As a result, the heterogeneity further promotes recovery: with the same perturbation size, only half of the domain recovers in the homogeneous initial setting, whereas the whole domain recovers when heterogeneity is present.

\subsection{On-branch final-state optimisation}
The third optimisation problem considered for the Allen-Cahn equation is the on-branch final-state optimisation. In contrast to the previous case, the system now starts from an untipped state. The scalar bifurcation parameter increases slowly over time with a prescribed linear time dependence, $\mu(t) = -1.5 + 5 \cdot 10^{-4} t$, across all simulations. With this gradual increase, the solution follows the stable branches of the bifurcation diagram throughout the evolution. 

As in the previous optimisation setting, the objective is to minimise the spatial mean of the solution which measures the spatial extent of the unfavourable state, at the final time $t_e$. This final time fixes the bifurcation parameter final value, denoted by $\mu_e = \mu(t_e)$. To achieve this, two constraints are imposed: first, the final state must result from the PDE evolution with time-varying $\mu(t)$, and second, the perturbation size is bounded, as before. Formally, we seek a perturbation $\mu_{\mathrm{pert}}$ that satisfies
\begin{align*}
    &\min_{\mu_{\mathrm{pert}}}  \langle y(x, t_e)\rangle_x \\
    &\text{subject to } 
    \begin{cases}
        -\dfrac{\partial y}{\partial t} + \dfrac{\partial^2 y}{\partial x^2} + y(1 - y^2) + \mu(t) + \mu_{\mathrm{het}}(x) + \mu_{\mathrm{pert}}(x) = 0,  \qquad t \in [t_0, t_e]\\[0.4em]
        y(x, t_0) = y_{\mathrm{IC}},\\[0.4em]
        \displaystyle \frac{1}{2}\int_{-1}^1 |\mu_{\mathrm{pert}}(x)|\,dx \le \delta.
    \end{cases}
\end{align*}
To illustrate this optimisation, we examine two scenarios that differ only in the final value of the bifurcation parameter: Case A takes $\mu_e = 0.63$, while Case B takes $\mu_e = 0.7$. In both cases, we fix the pre-existing heterogeneity to $\mu_{\mathrm{het}}(x) = \tfrac{1}{2} \cos(\pi x)$ and the perturbation size to $\delta = 0.3$. The results in \Cref{fig:3A} and \ref{fig:3B} show two distinct behaviours: in Case A, the optimisation prevents tipping altogether, while in Case B, it reduces the extent of tipping and leads to a partially fragmented final state.
\subsubsection{Case A: Complete tipping avoidance}
We first consider Case A, corresponding to a final bifurcation parameter value $\mu_e = 0.63$.
\begin{figure}
    \centering
    \includegraphics[width=\linewidth]{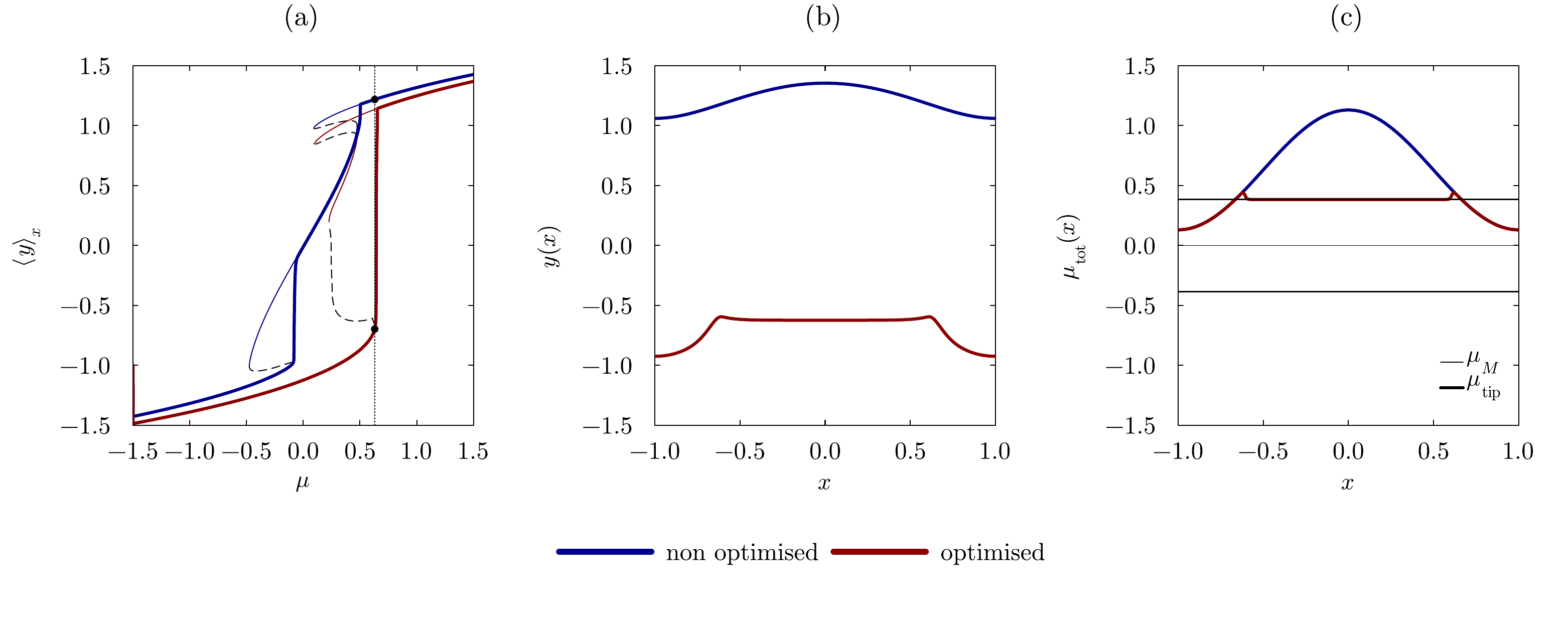}
    \caption{
On-branch final-state optimisation for $\mu_e = 0.63$ with heterogeneity $\mu_{\mathrm{het}}(x) = \tfrac{1}{2}\cos(\pi x)$ (Case A). 
\textbf{(a)} Bifurcation diagrams showing the spatial mean of the solution, $\langle y  \rangle_x = \frac{1}{2}\int_{-1}^{1}y(x)dx$, as a function of the bifurcation parameter $\mu$, plotted as thin curves (solid for stable branches, dashed for unstable branches). Superimposed thick curves show space-averaged, time-dependent trajectories obtained from numerical simulations for the non-optimised system (dark blue) and the optimised system (dark red). The vertical line indicates the value of the bifurcation parameter at the final simulation time, $\mu(t_{e})$.
\textbf{(b)} Spatial profiles $y(x)$ at the final simulation time, $t_e$, for the non-optimised (dark blue) and optimised (dark red) systems.
\textbf{(c)} Spatial heterogeneities $\mu_{\mathrm{tot}}(x)$ at $t_e$ for the non-optimised system (dark blue) and the optimised system (dark red). Horizontal lines indicate the local tipping points $\mu_{\mathrm{tip}}^{\pm}$ (thick line) and the Maxwell point $\mu_M$ (thin line).
}
    \label{fig:3A}
\end{figure}
\Cref{fig:3A} shows the results of Case A. In panel (a), the bifurcation diagrams and the space-averaged trajectories of the solutions are shown, with the optimised case in red and the non-optimised case in blue. The thin vertical line indicates $\mu_e$. The optimisation reduces the solution mean at $\mu_e$, from $1.2$ to $-0.7$. Panel (b) displays the solution at $\mu_e$ in both cases, plotted as a function of space. Panel (c) shows the total heterogeneities at $t=t_e$. In the non-optimised case (blue line), the solution tips in the central part of the domain as soon as the maximum of $\mu_{\mathrm{tot}}(x)$ crosses the homogeneous tipping point. As a consequence, two fronts are generated, and their locations match the intersection between $\mu_{\mathrm{tot}}(x)$ and the Maxwell point. Once $\mu_{\mathrm{tot}}(x)$ gets fully above the Maxwell point, all fronts disappear and the solution reaches the upper branch of the bifurcation diagram. The optimisation outputs a perturbation that keeps the total heterogeneity $\mu_{\mathrm{tot}}(x)$ below the tipping point until $\mu$ reaches $\mu_e$, thereby fully avoiding tipping.

The value of the scalar bifurcation parameter $\mu_e = 0.63$ corresponds to the output of the equilibrium optimisation problem with the same heterogeneity and perturbation size (\Cref{fig:1A}). Therefore, in this setting, tipping can be prevented entirely for any $\mu_e$ that is equal to or below the threshold obtained in the equilibrium optimisation, provided that the perturbation size $\delta$ is at least as large as that used to determine this threshold. By contrast, if $\mu_e$ is larger than $0.63$, some tipping is inevitable. To examine how this tipping can be mitigated, we now analyse Case B.

\subsubsection{Case B: Partial tipping mitigation beyond the equilibrium threshold}
We next consider Case B, with a larger final value $\mu_e = 0.7$, which lies beyond the threshold identified by equilibrium optimisation and therefore necessarily leads to some degree of tipping.
\begin{figure}
    \centering
    \includegraphics[width=\linewidth]{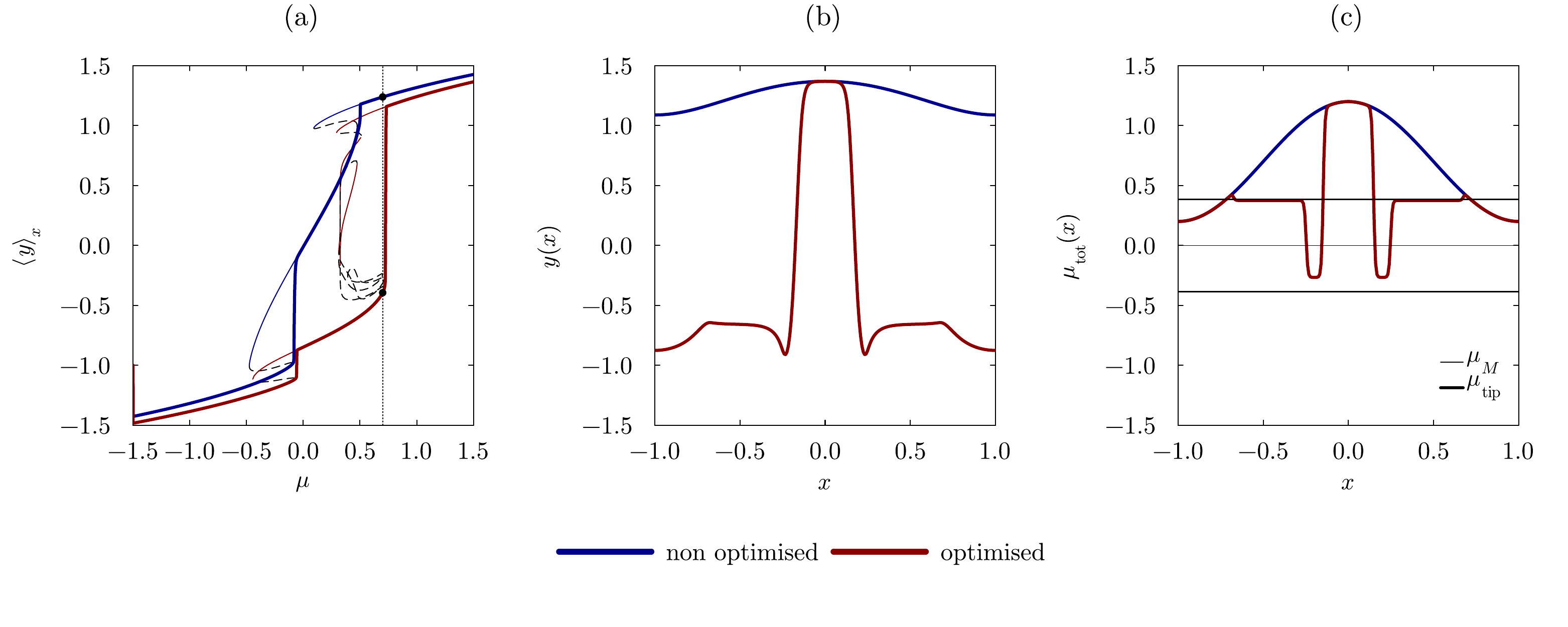}
    \caption{
On-branch final-state optimisation for $\mu_e = 0.7$ with heterogeneity $\mu_{\mathrm{het}}(x) = \tfrac{1}{2}\cos(\pi x)$ (Case B).
\textbf{(a)} Bifurcation diagrams showing the spatial mean of the solution, $\langle y  \rangle_x = \frac{1}{2}\int_{-1}^{1}y(x)dx$, as a function of the bifurcation parameter $\mu$, plotted as thin curves (solid for stable branches, dashed for unstable branches). Superimposed thick curves show space-averaged, time-dependent trajectories obtained from numerical simulations for the non-optimised system (dark blue) and the optimised system (dark red). The vertical line indicates the value of the bifurcation parameter at the final simulation time, $\mu(t_{e})$.
\textbf{(b)} Spatial profiles $y(x)$ at the final simulation time, $t_e$, for the non-optimised (dark blue) and optimised (dark red) systems.
\textbf{(c)}  Spatial heterogeneities $\mu_{\mathrm{tot}}(x)$ at $t_e$ for the non-optimised system (dark blue) and the optimised system (dark red). Horizontal lines indicate the local tipping points $\mu_{\mathrm{tip}}^{\pm}$ (thick line) and the Maxwell point $\mu_M$ (thin line).
}
    \label{fig:3B}
\end{figure}
\Cref{fig:3B} shows the results for $\mu_e = 0.7$. As discussed above, the optimisation cannot avoid tipping in this case, but it can limit its spatial extent. In panel (a), the optimised solution (red plot) shows a small tipping event at $\mu = -0.12$, followed by a second one at $\mu = \mu_e$. Panel (b) shows that the optimised solution at $\mu_e$ is only partially tipped: the alternate state appears only in the central region of the domain, contrary to the non-optimised case, where the solution is tipped everywhere.

To better understand this outcome, panel (c) displays the total heterogeneity. In the central region, the perturbation does not modify the pre-existing heterogeneity. Since the total heterogeneity exceeds $\mu_{\mathrm{tip}}^+$ there, the solution tips locally, forming two fronts. To prevent the expansion of these fronts, the perturbation lowers heterogeneity below the Maxwell point at the front locations, thereby preventing their propagation. In the outer region, the perturbation is applied only where the pre-existing heterogeneity is above the homogeneous tipping point, ensuring that the total heterogeneity remains below $\mu_{\mathrm{tip}}^+$ and no further tipping occurs. Notably, regions where the total heterogeneity lies between $\mu_M$ and $\mu_{\mathrm{tip}}^+$ tip locally in the non-optimised case, but remain untipped under optimisation, even though no perturbation is directly applied there.

\section{Discussion}\label{sec:5}
In this paper, we have proposed a general optimisation framework for spatially extended multistable systems. The framework is formulated as a constrained mathematical optimisation problem that identifies spatially targeted interventions to improve resilience under explicit dynamical and size constraints. Throughout this work, resilience is used operationally and quantified through setting-dependent objective functions.
The framework is illustrated using the one-dimensional Allen-Cahn equation \cite{allen1972_allen_cahn} with pre-existing heterogeneity as in \eqref{eq:AC_general}. This minimal conceptual model captures key mechanisms of gradient systems, including tipping and front pinning, and serves as a simple prototype for tipping in spatially extended systems \cite{bastiaansen2022_fragmented_tipping}.

Our results show that targeted small-scale spatial interventions can reshape larger-scale system behaviour by modifying both tipping thresholds and front dynamics \cite{krakovska2024_resilience_dynamical_systems, bastiaansen2020_effect_climchange_resilience}, often by localised interventions. We analysed three optimisation settings: one equilibrium optimisation and two final-state optimisations. The equilibrium optimisation shifts the saddle nodes to higher parameter values, extending the range of external pressure for which the favourable state exists, thus providing means to postpone tipping. The off-branch final-state optimisation creates, pins and reverses fronts, thereby promoting recovery of the favourable state or limiting the spread of the undesirable state. The on-branch final-state optimisation confines collapse to limited regions when tipping is unavoidable, preserving part of the domain in the desirable state. Together, these results illustrate how localised interventions can reshape global system behaviour.

The optimisation results can be understood through two core mechanisms: local crossings of tipping thresholds and of the Maxwell point by the total heterogeneity $\mu_{\mathrm{tot}}(x)$. These crossings determine whether local tipping occurs and whether fronts advance, remain pinned, or retreat. By modifying $\mu_{\mathrm{tot}}(x)$  where these critical thresholds are crossed or nearly crossed, the intervention can prevent tipping and control front dynamics. Pre-existing heterogeneity influences resilience through the same mechanisms, in particular by creating locations where fronts naturally pin or recovery can spread. More generally, our optimisation framework highlights a key principle for intervention in spatially extended tipping systems. Rather than requiring widespread modifications across the domain, resilience can be improved most effectively by targeting specific spatial regions. This reflects the inherently non-local nature of the dynamics: small perturbations applied over narrow intervals can determine the fate of the entire domain.

A central feature of our framework is its flexibility. In this paper, we chose a relatively simple model to clearly illustrate the underlying mechanisms, but the framework accommodates a wide range of extensions. First, the framework supports different intervention types: while current work uses additive perturbations, non-additive perturbations can also be implemented. 
Second, the objective function can be adapted to encode system-specific notions of resilience. For instance, in final-state optimisations, we use the global mean of the final state, but alternatives include local objectives or alternative resilience measures.
Third, constraints on the intervention can be tailored to reflect practical limitations of the system. For instance, spatial weighting can be used to encode spatially varying costs or feasibility, making interventions more restricted or expensive in certain regions.
Additional examples illustrating these extensions are provided in the supplementary material.
Finally, our framework supports pre-existing spatial heterogeneity and time-dependent parameters, making it applicable to a broad class of spatially extended systems.

Future work could develop the framework and its applications in several directions. One priority is to more explicitly quantify resilience \cite{krakovska2024_resilience_dynamical_systems} after optimisation. A resilience metric would measure how well the system withstands perturbations once an intervention is applied. A second direction is to apply the framework to non-gradient systems, since most realistic climate and ecological models lack an energy-like functional. In these models, the global potential does not exist, and, as a consequence, the Maxwell point may be absent or less clearly defined \cite{banerjee2026_rethinking_spatial_tipping}; fronts may also display more complex dynamics \cite{banerjee2026_rethinking_spatial_tipping, carter2023_2d_front_destabilization, FernandezOto2019_reverse_desertification}. A further extension is to apply the framework to two-dimensional pattern-forming systems \cite{rietkerk2002_selforganized_vegetation, meron2016_pattern_missing_link, meron2018_dryland_patterns, carter2023_2d_front_destabilization, FernandezOto2019_reverse_desertification}, which may display stripes, gaps, spots, and front instabilities such as fingering. In these settings, the framework could help design interventions that exploit front curvature or steer the system towards desirable spatial patterns. Finally, the framework could be extended to include time-dependent perturbations. Allowing the optimal intervention to vary over time would enable pulsed actions for short-term restoration or moving interventions along critical regions. It would naturally connect the static optimisation approach to an optimal control approach \cite{troltzsch2010_optimal, CasasRyllTroltzsch2013_opt_control}.

Looking ahead, the framework could be developed into a practical tool for identifying spatial strategies that reduce the risk of dangerous tipping events in ecosystems and climate subsystems. By combining data-informed spatial inputs, richer process-based models, and realistic representations of costs and feasibility, it could help evaluate where limited interventions would be most effective for strengthening resilience and promoting recovery under heterogeneous conditions.

\section*{Data availability statement}
No empirical datasets were created or analysed in this study. 
The Julia code used to generate the numerical simulations, optimisation results, and figures is available at \url{https://github.com/AuroraFaureRagani/optimisation_spatially_extended_systems}.
\section*{Acknowledgements}
This publication is part of the project ‘Resilience optimisation of spatially extended tipping prone systems’ with file number ‘OCENW.M.22.471’ of the research programme ‘Open Competition ENW’ which is (partly) financed by the Dutch Research Council (NWO).
\\
The authors declare no competing interests.
\printbibliography

\newpage
\appendix

\section{Additional optimisation variants}
\textbf{Description}: Additional numerical examples illustrating the flexibility of the proposed optimisation framework.

This supplementary material examines the flexibility of the optimisation framework beyond the specific formulations considered in the main text. It presents three additional numerical variants that modify, respectively, the objective function (\Cref{SMsec:local obj}), the intervention constraint (\Cref{SMsec: weighted cost}), and the mechanism through which the intervention enters the governing equation (\Cref{SMsec: non add pert}). While the applications in the main text use additive perturbations together with either equilibrium or global final-state objectives, the examples below show that the framework can accommodate alternative formulations at each of these three levels. These examples are intended as flexibility checks rather than as separate mechanistic analyses. Unless stated otherwise, the model, perturbation parametrisation, optimisation algorithm, budget value, and final time are the same as in Case B of the on-branch final-state optimisation presented in section 4.3.2 and figure 9 of the main text.

\subsection{Local objective function}\label{SMsec:local obj}

To illustrate that the framework is not restricted to global resilience measures, we repeat the on-branch final-state optimisation of section 4.3 using a local objective function. For a prescribed target region \(\Omega_{\mathrm{obj}}\subset[-1,1]\), we minimise
\[
F_{\Omega_{\mathrm{obj}}}(\mu_{\mathrm{pert}}) = 
\frac{1}{|\Omega_{\mathrm{obj}}|}
\int_{\Omega_{\mathrm{obj}}} y(x,t_e)\,dx .
\]
Since the desirable state corresponds to the lower branch, minimising this
quantity prioritises preservation or recovery of the desirable state inside
\(\Omega_{\mathrm{obj}}\), rather than across the full domain. All other components of the optimisation problem, including the perturbation parametrisation and the intervention budget, are kept unchanged.

\begin{figure}
    \centering
    \includegraphics[width=\linewidth]{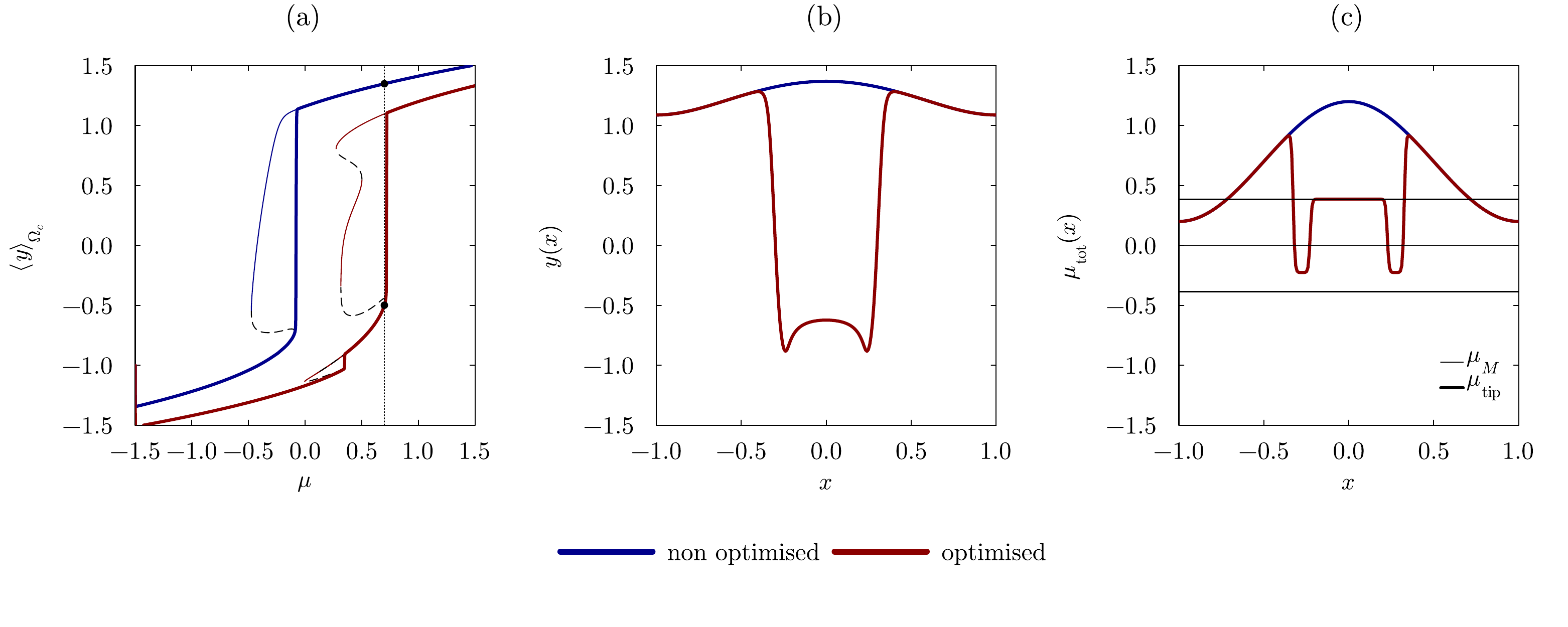}
    \caption{Optimisation results for the local objective function case.
    \textbf{(a)} Bifurcation diagrams showing the spatial mean of the solution over \(\Omega_{\mathrm{obj}}=[-1/3,1/3]\), $\langle y  \rangle_{\Omega_{\mathrm{obj}}} = \frac{3}{2}\int_{-1/3}^{1/3}y(x)dx$, as a function of the bifurcation parameter $\mu$, plotted as thin curves (solid for stable branches, dashed for unstable branches). Superimposed thick curves show space-averaged, time-dependent trajectories obtained from numerical simulations for the non-optimised system (dark blue) and the optimised system (dark red). The vertical line indicates the value of the bifurcation parameter at the final simulation time, $\mu(t_{e})$.
    \textbf{(b)} Spatial profiles $y(x)$ at the final simulation time, $t_e$, for the non-optimised (dark blue) and optimised (dark red) systems.
    \textbf{(c)}  Spatial heterogeneities $\mu_{\mathrm{tot}}(x)$ at $t_e$ for the non-optimised system (dark blue) and the optimised system (dark red). Horizontal lines indicate the local tipping points $\mu_{\mathrm{tip}}^{\pm}$ (thick line) and the Maxwell point $\mu_M$ (thin line).}\label{SMfig:loc obj}
\end{figure}

An example with \(\Omega_{\mathrm{obj}}=[-1/3,1/3]\), a central subdomain, is shown in \Cref{SMfig:loc obj}. Because the objective is evaluated only inside \(\Omega_{\mathrm{obj}}\), the optimisation sacrifices some global improvement in order to optimise the state inside the prescribed target region. Consequently, the optimal intervention favours front pinning near the boundaries of \(\Omega_{\mathrm{obj}}\). This illustrates how the same optimisation framework can encode spatially local notions of resilience.

\subsection{Spatially weighted intervention cost}\label{SMsec: weighted cost}

We next modify the intervention constraint used in the on-branch final-state optimisation in section 4.3 of the main text by introducing a spatially heterogeneous cost. Instead of penalising all perturbations equally across the domain, we impose the weighted budget constraint
\[
\frac{1}{2}\int_{-1}^{1} w(x)\,|\mu_{\mathrm{pert}}(x)|\,dx
\leq \delta,
\qquad w(x)\geq 0 .
\]
The weight \(w(x)\) represents the relative cost of applying an intervention
at location \(x\), which may encode differences in accessibility or
feasibility. Larger values of \(w\) make interventions more expensive in that part of the domain. To keep the weighted and unweighted budgets comparable, we normalise the weight such that
\[
\frac{1}{2}\int_{-1}^{1} w(x)\,dx = 1,
\]
so that the spatial average of \(w\) over the domain is one.

\begin{figure}
    \centering
    \includegraphics[width=\linewidth]{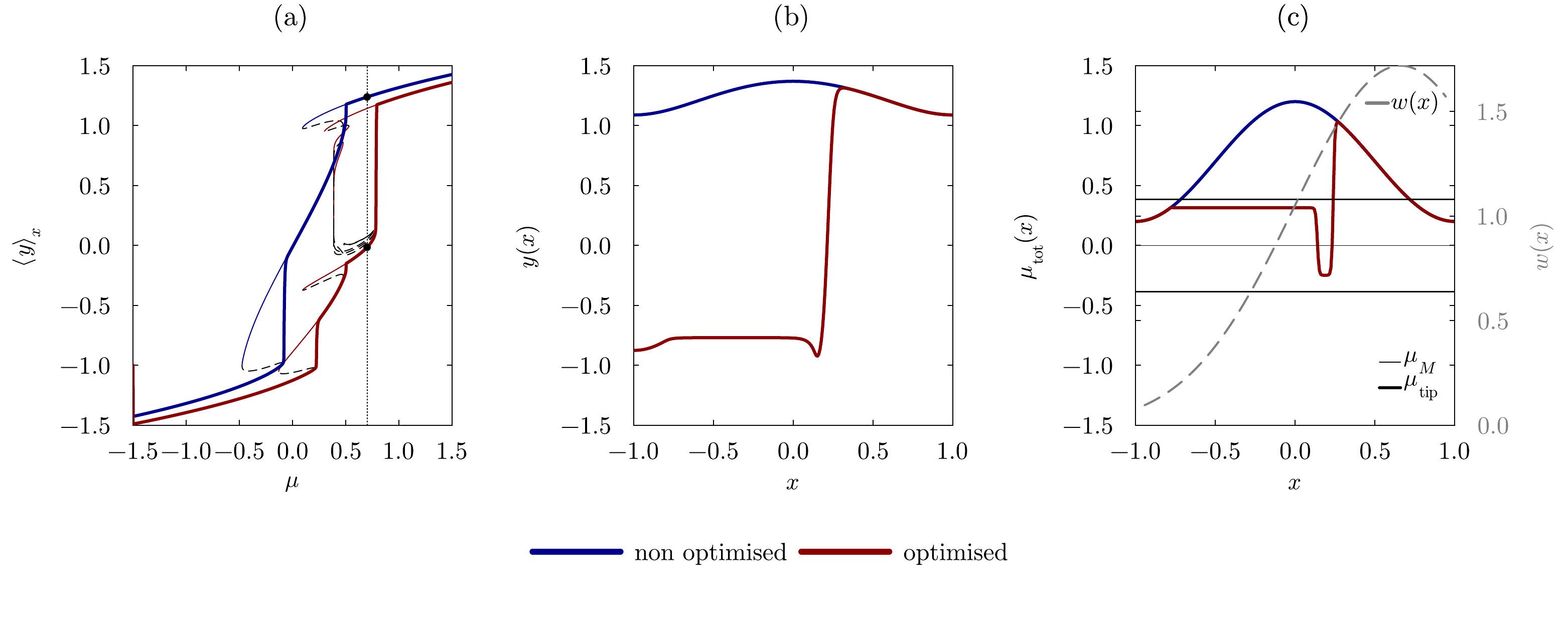}
    \caption{Optimisation results for the spatially weighted cost case.
    \textbf{(a)} Bifurcation diagrams showing the spatial mean of the solution over $[-1,1]$, $\langle y  \rangle_x = \frac{1}{2}\int_{-1}^{1}y(x)dx$, as a function of the bifurcation parameter $\mu$, plotted as thin curves (solid for stable branches, dashed for unstable branches). Superimposed thick curves show space-averaged, time-dependent trajectories obtained from numerical simulations for the non-optimised system (dark blue) and the optimised system (dark red). The vertical line indicates the value of the bifurcation parameter at the final simulation time, $\mu(t_{e})$.
    \textbf{(b)} Spatial profiles $y(x)$ at the final simulation time, $t_e$, for the non-optimised (dark blue) and optimised (dark red) systems.
    \textbf{(c)}  Spatial heterogeneities $\mu_{\mathrm{tot}}(x)$ at $t_e$ for the non-optimised system (dark blue) and the optimised system (dark red). Horizontal lines indicate the local tipping points $\mu_{\mathrm{tip}}^{\pm}$ (thick line) and the Maxwell point $\mu_M$ (thin line). The dashed grey curve shows the spatial weight $w(x)$, plotted against the secondary vertical axis on the right.}
    \label{fig:appendix_variants}
\end{figure}

The resulting optimisation is shown in \Cref{fig:appendix_variants}, with $w(x) = C \cdot e^{-(x-0.7)^2}$, where $C$ is the normalisation constant. The optimisation adapts the placement of the perturbation to the spatial cost landscape, avoiding expensive regions when possible and reallocating the intervention towards cheaper locations. Nevertheless, the qualitative dynamical mechanism remains the same as in the main text: the intervention acts by shifting the final forcing relative to local tipping thresholds and by controlling front propagation and pinning. This example shows that practical considerations, such as spatially varying accessibility or feasibility, can be encoded through a heterogeneous cost function without changing the underlying optimisation framework.

\subsection{Non-additive interventions}\label{SMsec: non add pert}
The previous two variants modified the objective function and the intervention constraint while keeping the additive intervention mechanism used in the main text. We now modify the mechanism through which the intervention enters the governing equation.
In the main text, interventions were introduced as additive perturbations of the local bifurcation parameter, so that the total forcing was given by 
\[
\mu_{\mathrm{tot}}(x)
=
\mu(t)+\mu_{\mathrm{het}}(x)+\mu_{\mathrm{pert}}(x).
\]
Here, instead, we consider the equation
\[
y_t
=
D y_{xx}
+
(y+z(x))(1-y^2)
+
\mu
+
\mu_{\mathrm{het}}(x).
\]
In this formulation, the intervention $z(x)$ changes the nonlinear dependence of the reaction term on the state $y$, rather than shifting the forcing $\mu+\mu_{\mathrm{het}}(x)$ additively. It is therefore no longer represented by $\mu_{\mathrm{pert}}(x)$.

To interpret the effect of the intervention, we define the local forcing as
\[
\eta(x, t)=\mu(t)+\mu_{\mathrm{het}}(x).
\]
At a fixed spatial location, and neglecting diffusion, the local dynamics therefore depend on the forcing $\eta$ and on the local intervention $z = z(x)$. For each fixed value of $z$, the upper local tipping threshold, at which the lower stable equilibrium disappears, is
\[
\eta_{\mathrm{tip}}(z)
:= \frac{2}{27}\left((z^2 + 3)^{3/2} - z(z^2 - 9)\right).\]
recovering the upper tipping threshold of the original Allen–Cahn equation introduced in section 4 of the main text.
The corresponding local Maxwell threshold, at which the lower and upper stable equilibria have equal potential, is
\[
\eta_M(z)
:=
-\frac{2}{3}z+\frac{2}{27}z^3 ,
\]
For $z = 0$, this gives $\eta_M(0) = 0$, recovering the Maxwell point discussed in section 2 of the main text.

Returning to the spatial problem, the local forcing $\eta(x,t)$ is compared with the thresholds evaluated at the local intervention z(x). Local loss of the lower state occurs when $\eta(x, t) > \eta_{\mathrm{tip}}(z(x))$,
whereas locations satisfying $\eta(x, t) \approx\eta_M(x)$ are candidate front-pinning locations.
The intervention therefore acts by changing the local tipping and Maxwell thresholds against which the forcing is compared; it does not change $\eta(x,t)$ itself. Equivalently, these conditions can be expressed in terms of the global bifurcation parameter $\mu$ by defining the spatially dependent thresholds
\[
\mu_{\mathrm{tip}}(x)
=
\eta_{\mathrm{tip}}(z(x))-\mu_{\mathrm{het}}(x),
\qquad \text{and} \qquad
\mu_M(x)
=
\eta_M(z(x))-\mu_{\mathrm{het}}(x).
\]
Hence, \(\mu(t)>\mu_{\mathrm{tip}}(x)\) indicates local loss of the lower state in the local problem, while \(\mu(t)\approx\mu_M(x)\) identifies candidate front-pinning locations. These thresholds are used only to interpret the local mechanisms underlying the optimised solution. The optimisation itself is performed using the full PDE and therefore retains diffusion, boundary effects, and interactions between fronts.

\begin{figure}
    \centering
    \includegraphics[width=\linewidth]{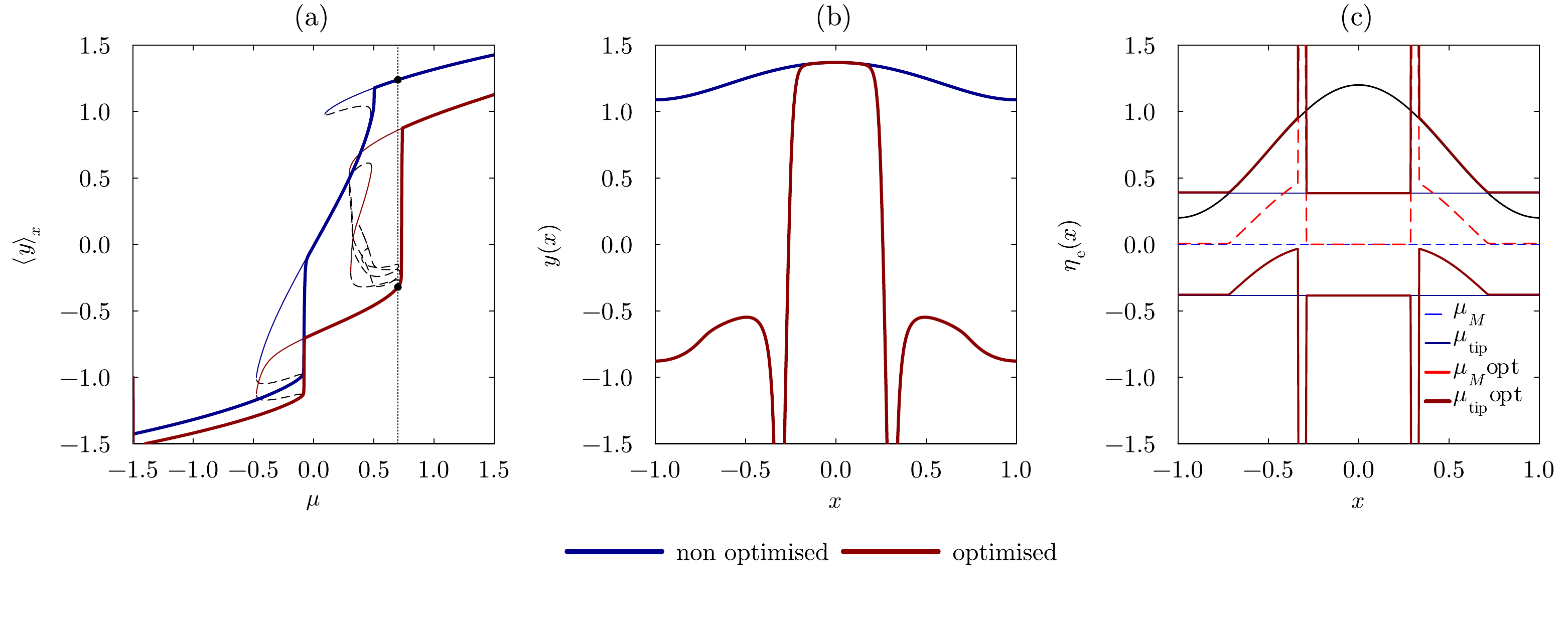}
    
    \caption{
    Final-state optimisation for the Allen--Cahn equation with non-additive intervention and heterogeneity $\mu_{\mathrm{het}}(x) = \tfrac{1}{2}\cos(\pi x)$.
    \textbf{(a)} Bifurcation diagrams showing the spatial mean of the solution over $[-1,1]$, $\langle y  \rangle_x = \frac{1}{2}\int_{-1}^{1}y(x)dx$, as a function of the bifurcation parameter $\mu$, plotted as thin curves (solid for stable branches, dashed for unstable branches). Superimposed thick curves show space-averaged, time-dependent trajectories obtained from numerical simulations for the non-optimised system (dark blue) and the optimised system (dark red). The vertical line indicates the value of the bifurcation parameter at the final simulation time, $\mu(t_{e})$.
    \textbf{(b)} Spatial profiles $y(x)$ corresponding to the parameter value $\mu = \mu_{e}$ for the non-optimised (dark blue) and optimised (dark red) systems.
    \textbf{(c)} Final local forcing \(\eta_{e}(x)\), together with the relevant local tipping and Maxwell thresholds. The dark blue lines show the thresholds for the non-optimised system with \(z=0\), while the dark red curves show the spatially dependent thresholds induced by the optimised non-additive intervention \(z(x)\). }\label{fig:nonadditive_intervention}
\end{figure}

The intervention parametrisation used in \Cref{fig:nonadditive_intervention} is constructed around these two mechanisms: postponing local loss of the desirable state and controlling front propagation.
At the final time $t_e$, we define $\mu_e = \mu(t_e)$ and $\eta_e(x) = \eta(x,t_e) = \mu_e +\mu_{\mathrm{het}}(x)$.
In an unconstrained setting, one could prevent local loss of the lower state by choosing \(z(x)\) so that the local tipping threshold matches the final forcing,
\begin{equation}\label{eq:fold_match}
\eta_{\mathrm{tip}}(z(x)) = \eta_{e}(x).
\end{equation}
This condition places the local saddle-node at the final forcing, so that the lower equilibrium remains locally available for \(\mu(t) < \mu_{e}\).

Under a finite intervention budget, however, imposing this condition throughout the domain may be impossible or suboptimal. No intervention is required where the lower state is already locally available without intervention, namely where
\[
\eta_{e}(x)\leq \eta_{\mathrm{tip}}(0).
\]
If the remaining budget is still insufficient to protect all locally vulnerable regions, the optimisation may instead allow part of the domain to tip and use the intervention $z(x)$ to control the fronts generated by this local transition. In this case, the intervention not only shifts the local tipping threshold but also adjusts the local Maxwell threshold so that
\[
\eta_M(z(x))\approx \eta_{e}(x),
\]
at selected locations, thereby favouring front pinning.

For the cosine heterogeneity used in Case B of the on-branch final-state optimisation in section 4.3.2 and figure 9 of the main text
\[
\mu_{\mathrm{het}}(x)=\frac{1}{2}\cos(\pi x),
\]
the final local forcing \(\eta_{e}(x)\) is largest in the central part of the domain. Under the finite budget considered here, the optimisation leaves a central region insufficiently protected, thereby allowing it to tip. On either side of this region, the optimisation selects $z(x)$ so that the local Maxwell threshold $\eta_M(z(x))$ is brought close to the fixed local forcing $\eta_e(x)$, thereby arresting the resulting fronts. 

In the outer parts of the domain, the intervention shifts the local tipping threshold towards the final forcing, preventing independent local loss of the lower state away from the fronts. The optimised degrees of freedom therefore consist of the size of the central tipping region, the magnitude of $z(x)$ in the front-pinning regions, and the extent of the outer regions in which the fold-matching condition \eqref{eq:fold_match} is safistied.

The resulting optimisation is shown in figure \Cref{fig:nonadditive_intervention}. This example changes the intervention mechanism while leaving the overall optimisation procedure unchanged. Once z(x) has been parametrised, the optimisation proceeds as in the additive case, and the resulting solution remains interpretable in terms of local tipping and Maxwell thresholds.

Together, these variants show that the optimisation framework is flexible at
three distinct levels: the objective functional, the intervention constraint,
and the intervention mechanism itself. In each case, the numerical optimisation
is performed within the same PDE-constrained framework, while the resulting
solutions remain interpretable in terms of local thresholds and Maxwell points.

\end{document}